\documentclass[12pt]{article}

\author{Mathew D.~Penrose\footnote{e-mail: \texttt{m.d.penrose@bath.ac.uk}}
and Andrew R.~Wade\footnote{e-mail: \texttt{a.wade@bath.ac.uk}}\\
\normalsize
 Department of Mathematical Sciences,
 University of Bath,\\
\normalsize
 Bath BA2 7AY, England.
}
\title{Limit theory for the random on-line nearest-neighbour graph}

\date{March 2006}

\usepackage{amssymb}
\usepackage{epsfig}
\usepackage{color}
\newcommand{\ud}{\mathrm{d}}
\newtheorem{theorem}{Theorem}[section]
\newtheorem{proposition}{Proposition}[section]
\newtheorem{lemma}{Lemma}[section]

\newcommand{\bean}{\begin{eqnarray*}}
\newcommand{\eean}{\end{eqnarray*}}
\newcommand{\bea}{\begin{eqnarray}}
\newcommand{\eea}{\end{eqnarray}}
\newcommand{\tod}{\stackrel{{\cal D}}{\longrightarrow}}
\newcommand{\eqd}{\stackrel{{\cal D}}{=}}

\newcommand{\toas}{\stackrel{{\rm a.s.}}{\longrightarrow}}
\newcommand{\inL}{\stackrel{L^1}{\longrightarrow}}
\newcommand{\inLL}{\stackrel{L^2}{\longrightarrow}}
\newcommand{\inLLL}{\stackrel{L^3}{\longrightarrow}}

\newcommand{\rem}{\noindent \textbf{Remark. }}
\newcommand{\rems}{\noindent \textbf{Remarks. }}
\newcommand{\proof}{\noindent \textbf{Proof. }}
\def\Exp{E}
\def\Pr{P}
\def\Var{{\mathrm{Var}}}
\def\Cov{{\rm Cov}}
\def\R{{\bf R }}

\def\1{{\bf 1 }}
\def\N{{\bf N }}

\def\X{{\cal X}}

\def\U{{\cal U}}
\def\TT{{\cal T}}
\def\OO{{\cal O}}
\def\tO{\tilde {\cal O}}
\def\Po{{\cal P}}
\def\H{{\cal H}}
\def\Jal{J_\alpha}
\def\tJal{\tilde J_\alpha}
\def\Hal{H_\alpha}
\def\tHal{\tilde H_\alpha}
\def\tGal{\tilde G_\alpha}

\def\tJo{{\tilde J_1}}
\def\tHo{{\tilde H_1}}
\def\tGo{{\tilde G_1}}
\def\tY{{\tilde Y}}

\def\tR{{\tilde R}}
\def\tS{{\tilde S}}
\def\LL{{\cal L}}

\def\onng{{\rm ONG}}

\def\bx{{\bf x}}
\def\by{{\bf y}}

\def\bU{{\bf U}}
\def\0{{\bf 0}}

\begin{document}

\maketitle
\abstract{ In the on-line nearest-neighbour graph (ONG), each point after the first in a
sequence of points in $\R^d$ is joined by an edge to its nearest-neighbour amongst 
those points that precede it in the sequence. 
We study the large-sample asymptotic behaviour of the 
total power-weighted length of the ONG on uniform random points in $(0,1)^d$. In
particular, for $d=1$ and weight exponent $\alpha>1/2$,
the limiting distribution of the centred
total weight is characterized by a distributional
fixed-point equation. As an ancillary result, we give exact expressions for the expectation and variance
of the standard nearest-neighbour (directed) graph on uniform random points in the unit interval.
}

\vskip 3mm

\noindent
{\em Key words and phrases:} Nearest neighbour graph; spatial network evolution;
weak convergence; fixed-point equation; divide-and-conquer.

\vskip 3mm

\noindent
{\em AMS 2000 Mathematics Subject Classification:} Primary: 60D05, 60F05; Secondary: 90B15.

\section{Introduction}

Spatial graphs, defined
on random point sets in Euclidean space, constructed
 by joining nearby points
according to some deterministic rule, have been the subject of considerable
recent interest. Examples of such graphs include the geometric graph,
the minimal-length spanning tree, and 
the nearest-neighbour graph and its relatives. Many aspects of the
large-sample asymptotic theory for such graphs, which are locally
determined in a certain sense, are by now quite well understood.
See for example
\cite{KL,penbook,mdp,penyuk1,penyuk2,steelebook,yukbook}.

Many real-world networks have several common features, including spatial structure,
local construction (nearby points are more likely to be connected), and sequential growth (the network
evolves over time via the addition of new nodes). 
In this paper our main object of interest is the {\em on-line nearest-neighbour graph}, 
which is one of the simplest models of network evolution that captures some of these features.
We give a detailed description later. Recently, graphs with an
`on-line' structure, i.e.~in which vertices are added
sequentially and connected to existing vertices via some rule,
have been the subject of considerable study in relation to the
modelling of real-world networks. The non-rigorous literature is
extensive (see for example \cite{dorog,newman} for surveys), but
rigorous mathematical results are fewer in number, even for simple
models, and the existing results concentrate on graph-theoretic
rather than geometric properties (see e.g.~\cite{bbcr,bol1}).

\begin{figure}[h!]
\begin{center}
\includegraphics[angle=0, width=0.5\textwidth]{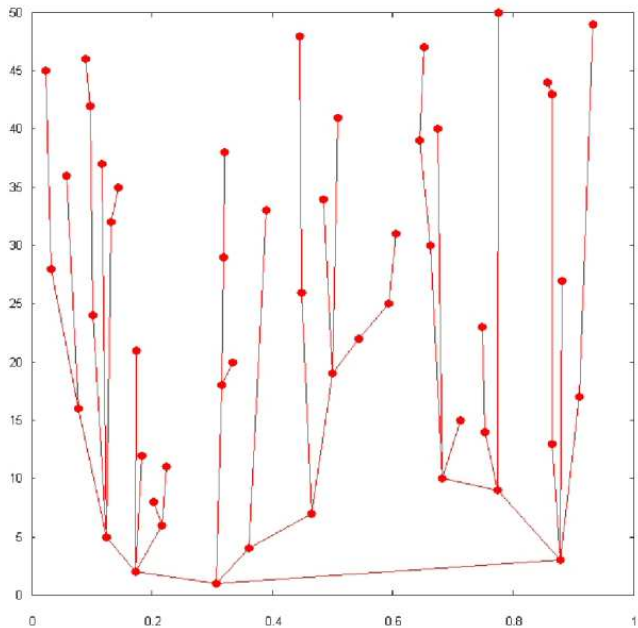}
\end{center}
\caption{Realization of the ONG
 on 50 simulated uniform random points in the unit interval. The vertical axis
gives the order in which the points arrive, and their position
is given by the horizontal axis.}
\label{ong1d}
\end{figure}

The on-line nearest-neighbour graph (or
$\onng$ for short) 
is constructed on $n$ points arriving sequentially in $\R^d$ by
connecting each point to its nearest neighbour amongst the
preceding points in the sequence.
The $\onng$ was apparently introduced
in \cite{bbcr} as a simple growth model
of the world wide web graph (for $d=2$). When $d=1$, the $\onng$
is related to certain fragmentation
processes, which are of separate interest in relation to, for example,
 molecular fragmentation (see
e.g.~\cite{bertoin}, and references therein). The $\onng$ in $d=1$
is related to the so-called `directed linear tree' considered in
\cite{total}.
The higher dimensional $\onng$ has also been studied \cite{mdp}. Figure \ref{ong1d}
shows a realization of the $\onng$ on $50$ simulated random points in the unit interval.
Figure \ref{onngfig} below shows realizations of the planar and three-dimensional
$\onng$, each on $50$ simulated uniform random points.

\begin{figure}[h!]
\begin{center}
\includegraphics[angle=0, width=0.9\textwidth]{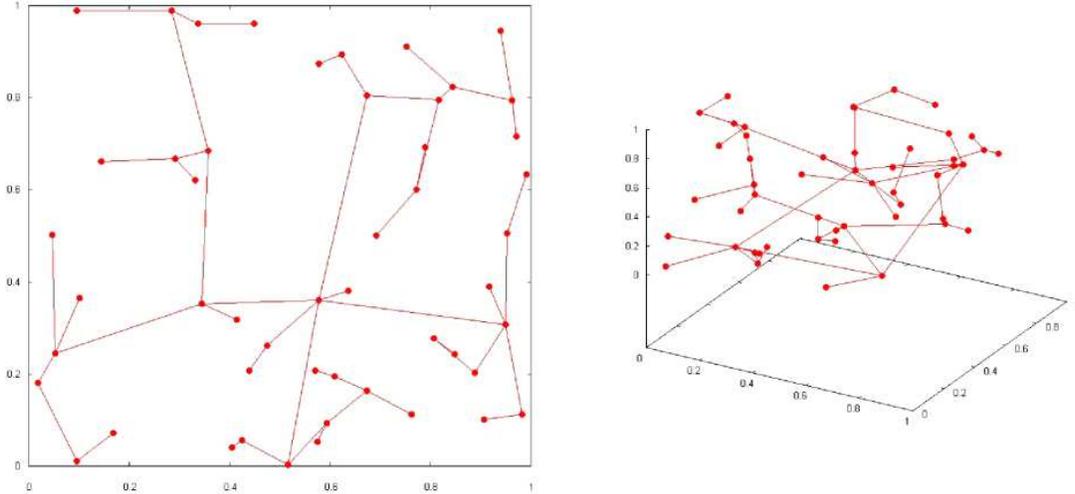}
\end{center}
\caption{Realizations of the ONG
 on 50 simulated uniform random points in the unit square (left)
and the unit cube (right).}
\label{onngfig}
\end{figure}

We consider the total power-weighted length of the ONG on uniform random points
in $(0,1)^d$, $d \in \N$. We are interested in large-sample asymptotics, as
the number of points tends to infinity. Explicit laws of large numbers for the random ONG in
$(0,1)^d$ are given in \cite{llnpaper}. In the present
paper we give
further results on the limiting behaviour in general dimensions $d$. 

The main part of the present paper is concerned with convergence in distribution results for the
ONG. We give detailed properties
of the random ONG on uniform random points
in the unit interval ($d=1$), and identify the limiting
distribution of the centred total power-weighted length of the graph. When the
weight exponent $\alpha$ is greater than $1/2$, this distribution
is described in terms of a
 distributional fixed-point
equation reminiscent of those encountered in, for example, the
analysis of stochastic `divide-and-conquer' or recursive
algorithms. Such fixed-point distributional equalities, and the
recursive algorithms from which they arise, have received
considerable attention recently; see, for
example,~\cite{aldous,neinrusch,rosler0,rosler}. 

On the other
hand, we believe that for $\alpha \in (0,1/2]$ the total weight,
suitably centred and scaled, satisfies a central limit theorem (CLT).
Penrose \cite{mdp} gave such a result
for $\alpha \in (0,1/4)$. We believe that it should be possible to derive the CLT
for all $\alpha \in (0,1/2]$ via the divide-and-conquer methods of this paper. 
The main difficulty is to show that the variance of the total weight of the graph
scales appropriately in the large sample limit. We hope to address this in future work.

In this paper we also give new explicit results on the expectation and variance
of the standard one-dimensional nearest-neighbour (directed)
graph, in which each point is joined by a directed edge to its
nearest-neighbour, on uniform random points in the unit interval. This is related to our results on the
one-dimensional $\onng$ via the theory of Dirichlet spacings, which we make use of in our analysis. 

\section{Definitions and main results} \label{con}

Let $\X$ be a finite sequence of points in $\R^d$, and let $\| \cdot\|$
be the Euclidean norm. 
For $d \in \N$, let
\bea
\label{0818c}
v_d
: = \pi^{d/2} \left[ \Gamma \left( 1+ (d/2) \right) \right]^{-1},
\eea
 the volume of the unit
$d$-ball (see e.g.~equation (6.50) of \cite{huang}).

Define $w$ to be a weight function
on edges, assigning weight
$w(\bx,\by)$ to the edge
between $\bx\in \R^d$ and $\by \in \R^d$, such that $w: \R^d \times \R^d \to [0, \infty)$.
A case of particular interest is when the weight is taken to be
power-weighted Euclidean distance. In this case, for some $\alpha \geq 0$,
we have the weight function
\bea
\label{wf}
w_\alpha (\bx,\by) := \| \bx-\by \|^\alpha ,\eea
for $\bx, \by \in \R^d$.

\subsection{The on-line nearest-neighbour graph}
\label{subseconng}

We now give a formal definition of the on-line nearest-neighbour graph (ONG). 
Let $d \in \N$. Suppose $\bx_1, \bx_2, \ldots$
are points in $(0,1)^d$, arriving
sequentially;
the ONG on vertex set $\{ \bx_1,\ldots,\bx_n\}$
is formed by connecting each
point $\bx_i$, $i=2,3,\ldots,n$ to its nearest neighbour (in the Euclidean sense)
amongst
the preceding points in the sequence (i.e.~$\bx_1, \ldots, \bx_{i-1}$),
using the lexicographic ordering on $\R^d$ to break any ties. We call the
 resulting
tree
the ONG on $(\bx_1,\bx_2,\ldots,\bx_n)$.

From now on we take the sequence of points to be {\em
random}. Let $\bU_1, \bU_2,\ldots$ be a sequence of independent
uniform random vectors on $(0,1)^d$. Then for $n \in \N$ take
$\U_n=(\bU_1,\bU_2,\ldots,\bU_n)$, the binomial
point process consisting
of $n$ independent uniform random vectors
on $(0,1)^d$. Denote the ONG constructed on $\U_n$
by $\onng(\U_n)$. We restrict our analysis to the case of uniformly
distributed points.
Note that, with
probability one, $\U_n$
has
distinct
inter-point distances
so that the ONG on $\U_n$ is almost
surely unique.

The $\onng$ is of interest as a natural growth model for random
spatial graphs; in particular it has been used (with $d=2$) in the
context of the world wide web graph (see \cite{bbcr}). In
\cite{mdp}, stabilization techniques were used to prove that the
total length (suitably centred and scaled) of the $\onng$ on
uniform random points in $(0,1)^d$ for $d > 4$ converges in
distribution to a normal random variable. It is suspected that a
CLT also holds for $d=2,3,4$. On the other hand,
when $d=1$, the limit is not normal, as demonstrated by Theorem
\ref{onng1} (ii) below.

For $d \in \N$ and $\alpha \geq 0$, let $\OO^{d,\alpha} (\U_n)$
denote the total weight, with weight function $w_\alpha$ as given
by (\ref{wf}), of $\onng (\U_n)$. Our results for the $\onng$ in
general dimensions are as follows, and constitute a distributional convergence result for
$\alpha>d$, and asymptotic behaviour of the mean for $\alpha=d$. For the sake of
completeness, we include the 
law of large numbers for $\alpha<d$ from \cite{llnpaper} as part (i) of the theorem below.

\begin{theorem}
\label{onngthm}
Suppose $d \in \N$. We have the following:
\begin{itemize}
\item[(i)] Suppose
 $0 \leq \alpha <d$. Then, as $n \to \infty$
\bea
\label{0915az}
n^{(\alpha-d)/d} \OO^{d,\alpha} (\U_n)
\inL \frac{d}{d-\alpha} v_d^{-\alpha/d} 
\Gamma (1+(\alpha/d)).
\eea
\item[(ii)]
Suppose $\alpha>d$. Then, as $n \to \infty$,
\bea
\label{1027h}
\OO^{d,\alpha} (\U_n) \longrightarrow W(d,\alpha),
\eea
where the convergence is in $L^p$, $(p\in\N)$, and almost sure,
and $W(d,\alpha)$ is a nonnegative random variable with
$\Exp[(W(d,\alpha))^k]<\infty$ for $k\in\N$.
\item[(iii)]
Suppose $\alpha=d$. Then, as $n \to \infty$,
\bea
\label{1029a}
\Exp [\OO^{d,d} (\U_n)] = v_d^{-1} \log{n} +o(\log n).
\eea
\end{itemize}
\end{theorem}
In particular (\ref{1029a}) implies that
$\Exp[ \OO^{1,1} (\U_n)] \sim (1/2) \log n$, a result given more precisely
in Proposition \ref{propw0}
 below.
We prove Theorem \ref{onngthm} (ii) and (iii)
in Section \ref{prf1027}.

Now we consider the particular case of the $\onng$ in $d=1$, where $\U_n$
is now a sequence of independent uniform random points in the unit interval $(0,1)$.
 Let $\gamma$ denote Euler's constant,
so that $\gamma \approx 0.57721566$ and
\bea
\label{gamma}
\left( \sum_{i=1}^k \frac{1}{i} \right) - \log k = \gamma + O(k^{-1}).
\eea
The following result gives the expectation of the total weight of $\onng(\U_n)$.
\begin{proposition}
\label{propw0}
As $n \to \infty$, we have
\bean
 \Exp[ \OO^{1,\alpha}(\U_n) ] & = &
\frac{\Gamma(\alpha+1)}{1-\alpha} 2^{-\alpha}n^{1-\alpha} + \frac{2}{\alpha}-\frac{2^{-\alpha}(2-\alpha)}{\alpha(1-\alpha)}
+O(n^{-\alpha}); ~~~  (0<\alpha<1) \\
\Exp[ \OO^{1,1}(\U_n) ] & = & \frac{1}{2}
\log{n} + \frac{\gamma}{2}-\frac{1}{4} +o(1); ~~~  \\
\Exp[ \OO^{1,\alpha}(\U_n) ] & = & \frac{2}{\alpha (\alpha+1)} \left( 1 + \frac{2^{-\alpha}}{\alpha-1} \right)
+O(n^{1-\alpha}) ~~~ (\alpha>1)
 \eean
\end{proposition}
\proof The proposition follows from Proposition \ref{propw2} with Lemma \ref{0204a}. $\square$ \\

In Theorem \ref{onng1} below, we present our main convergence in
distribution results for the total
weight of the ONG (centred, in some cases) in $d=1$.
The limiting distributions
are of different types depending on the value of $\alpha$ in the weight
function (\ref{wf}). In this paper, we restrict attention to $\alpha>1/2$, and we
define
these limiting distributions in Theorem \ref{onng1}, in terms of distributional
fixed-point equations (sometimes called {\em recursive distributional
equations}, see \cite{aldous}). These fixed-point equations are of the form
\bea
\label{fixed}
X \eqd \sum_{r=1}^k A_r X^{\{r\}} + B,
\eea
where $k \in \N$, $X^{\{r\}}, r=1,\ldots,k$, are independent
copies of the random variable $X$, and $(A_1,\ldots,A_k,B)$
is a random vector, independent of $(X^{\{1\}},\ldots,X^{\{k\}})$,
satisfying the conditions
\bea
\label{fxdcnd}
\Exp \sum_{r=1}^k | A_r |^2 <1 , ~~~~~ \Exp [ B] = 0, ~~~~~
\Exp [ B^2 ] < \infty.
\eea
Theorem 3 of R\"osler \cite{rosler0} (proved by the
contraction mapping theorem; see also \cite{neinrusch,rosler})
says that if (\ref{fxdcnd}) holds, there is a unique
square-integrable distribution
with mean zero satisfying the fixed-point equation (\ref{fixed}),
and this will guarantee uniqueness of solutions to all the
distributional fixed-point equalities considered in the sequel.

We now define the distributions that will appear as limits
in Theorem \ref{onng1}, in terms of (unique) solutions
to fixed-point equations. In each case, $U$ denotes a uniform
random variable on $(0,1)$, independent of the other random variables
on the right hand side of the distributional equality. The fixed-point equations
(\ref{0919x})--(\ref{0923n}) are all of the form of (\ref{fixed}), and hence define unique
solutions.

We define $\tJo$ by the distributional fixed-point equation
\bea
\label{0919x}
\tJo \eqd \min \{ U, 1-U \} + U \tJo^{\{1\}} + (1-U) \tJo^{\{2\}}
+ \frac{U}{2} \log U + \frac{1-U}{2} \log (1-U) .
\eea
We shall see later (Proposition \ref{0201a}) that $\Exp [ \tJo] =0$.
For $\alpha>1/2$, $\alpha \neq 1$, define $\tJal$ by
\bea
\label{0919z}
\tJal \eqd U^\alpha \tJal^{\{1\}}
+ (1-U)^\alpha \tJal^{\{2\}} + \min \{ U^\alpha, (1-U)^\alpha \} + \frac{2^{-\alpha}}{\alpha-1}
\left( U^\alpha +(1-U)^\alpha -1 \right).
\eea
Define the random variable $\tHo$ by
\bea
\label{0924e}
\tHo \eqd U \tJo + (1-U) \tHo + \frac{U}{2} +
\frac{U}{2} \log U
+ \frac{1-U}{2} \log (1-U),
\eea
where $\tJo$ has the distribution given by (\ref{0919x}), and is independent
of the $\tHo$ on the right.
We shall see later (Theorem \ref{onng01}) that $\Exp [ \tHo] =0$. We give the first three
moments of $\tJo$ and $\tHo$ in Table \ref{tabmoms} later in this paper.
For $\alpha>1/2$, $\alpha \neq 1$, define $\tHal$ by
\bea
\label{0923n}
\tHal \eqd U^\alpha \tJal
+ (1-U)^\alpha \tHal + U^\alpha \left( 1+ \frac{2^{-\alpha}}{\alpha -1} \right)
+ ((1-U)^\alpha-1) \left( \frac{1}{\alpha} + \frac{2^{-\alpha}}{\alpha (\alpha-1)} \right),
\eea
where $\tJal$ has the distribution given
by (\ref{0919z}) and is independent of the $\tHal$ on the right.
We shall see later that, for $\alpha>1$,
the $\tJal$ and $\tHal$ defined in (\ref{0919z})
and (\ref{0923n})
arise as centred versions of the random variables $\Jal$ and $\Hal$, respectively,
satisfying the slightly simpler fixed-point equations (\ref{0919y}) and (\ref{0923l})
below, so that $\Exp [ \tJal] =\Exp[ \tHal]=0$; see
Proposition \ref{0920c}.
For $\alpha>1$, we have
\bea
\label{0919y}
\Jal \eqd U^\alpha \Jal^{\{1\}}
+ (1-U)^\alpha \Jal^{\{2\}} + \min \{ U^\alpha, (1-U)^\alpha \}.
\eea
Also for $\alpha>1$, we have
\bea
\label{0923l}
\Hal \eqd U^\alpha + U^\alpha \Jal + (1-U)^\alpha \Hal ,\eea
where $\Jal$ has distribution given by (\ref{0919y}) and is independent of the $\Hal$
on the right.
The expectations of $\Jal$ and $\Hal$ are given in Proposition \ref{0920c}. Note that the uniqueness
of the $\tJal$ and $\tHal$ implies the uniqueness of $\Jal$ and $\Hal$ also.

Theorem \ref{onng1} gives
our main results for the $\onng (\U_n)$ in one dimension. Theorem \ref{onng1}
will follow as a corollary to Theorem \ref{onng01}, which we present later.
Let $\tO^{d,\alpha} (\U_n) := \OO^{d,\alpha} (\U_n) - \Exp [\OO^{d,\alpha} (\U_n) ]$
be the centred total weight of the ONG on $\U_n$.
For ease of notation, we define the following random variables. As before,
$U$ is uniform on $(0,1)$ and independent of the other
variables on the right. For $1/2<\alpha<1$, let
\bea
\label{1006a}
 \tGal \eqd
U^\alpha \tHal^{\{1\}} \! + \!(1-U)^\alpha \tHal^{\{2\}} \!
 + \!
 \left( U^\alpha +(1-U)^\alpha - \frac{2}{1+\alpha} \right) \!\!
\left( \frac{1}{\alpha} - \frac{2^{-\alpha}}{\alpha(1-\alpha)} \right) \!\!,
\eea
where
$\tHal^{\{1\}}, \tHal^{\{2\}}$ are independent with distribution
given by (\ref{0923n}). Also let
\bea
\label{1006b}
 \tGo \eqd
U \tHo^{\{1\}}
+(1-U) \tHo^{\{2\}} + \frac{U}{2} \log U +\frac{1-U}{2}  \log (1-U) + \frac{1}{4},
\eea
where
$\tHo^{\{1\}}, \tHo^{\{2\}}$ are independent
with distribution given by (\ref{0924e}).  Now we state our
convergence in distribution results. We prove Theorem \ref{onng1}
in Section \ref{1d}.

\begin{theorem}
\label{onng1}
\begin{itemize}
\item[(i)]
For $1/2<\alpha<1$, we have that, as $n \to \infty$,
\bea
\label{1003a}
\tO^{1,\alpha} (\U_n)  \tod \tGal,
\eea
where $\tGal$ has distribution given by (\ref{1006a}), and $\Exp[\tGal]=0$.
\item[(ii)]
For $\alpha =1$, we have that, as $n \to \infty$,
\bea
\label{0202a}
\OO^{1,1} (\U_n) - \frac{1}{2} \left( \gamma + \log{n} \right) +\frac{1}{4} \tod
\tGo,
\eea
where $\tGo$ has distribution given by (\ref{1006b}). Also, $\Exp[\tGo]=0$,
$\Var [ \tGo]
= (19+4\log{2}-2\pi^2)/48 \approx 0.042362$,
and $\Exp[\tGo^3] \approx 0.00444287$.
\item[(iii)]
For $\alpha > 1$, the distribution
of the limit $W(1,\alpha)$ of (\ref{1027h})
is given by
\bean
W(1,\alpha)  \eqd
U^\alpha \Hal^{\{1\}}
+(1-U)^\alpha \Hal^{\{2\}},
\eean
where
$\Hal^{\{1\}}, \Hal^{\{2\}}$ are independent with the distribution
given by (\ref{0923l}).
\end{itemize}
\end{theorem}

\rems
(a) In Theorem 3.6 of \cite{mdp}, a CLT for
$\tO^{d,\alpha}(\U_n)$ is obtained for the case $0<\alpha<d/4$. In the context
of Theorem \ref{onngthm}, the result of \cite{mdp} implies that, provided $0<\alpha<d/4$,
as $n \to \infty$,
$n^{(\alpha/d)-(1/2)} \tO^{d,\alpha} (\U_n)$ is asymptotically normal.
In \cite{mdp}, it is remarked that it should be possible to extend the result to the case $d/4 \leq \alpha<d/2$
and perhaps $\alpha=d/2$ also. We hope
to address this in future work; in particular, the case $d=1$ should be amenable to solution via the divide-and-conquer
approach of this paper.

(b) A closely related `directed' version of the one-dimensional $\onng$ is the
`directed linear tree' (DLT) introduced in \cite{total}, in which each point is joined
to its nearest-neighbour to the {\em left} amongst those points preceding it in the sequence, if
such points exist. In \cite{total},
results for the DLT with $\alpha \geq 1$ analogous to parts (ii) and (iii) of Theorem \ref{onng1} were given.
Following the methods of the present paper, one can obtain results for the DLT with $1/2<\alpha<1$
analogous to part (i) of Theorem \ref{onng1}.

(c) Of interest is the limit behaviour of $\OO^{d,d} (\U_n)$ (i.e.~when $\alpha=d$).
When $d=1$, we have that $\OO^{1,1} (\U_n) - \Exp[\OO^{1,1} (\U_n)]$ converges in distribution
to a non-normal limiting random variable (see Theorem \ref{onng1} (ii)). It would be interesting to
determine whether $\OO^{d,d} (\U_n) - \Exp[\OO^{d,d} (\U_n)]$ converges in distribution
to a nondegenerate random variable for general $d=2,3,4,\ldots$, and whether or not this distribution is
normal.

(d) With some more detailed calculations (given in \cite{thesis}),
one can replace the error term $o(\log{n})$ in (\ref{1029a})
by $O(1)$ (see the remark in Section \ref{prf1027}).

\begin{figure}[h!]
\begin{center}
\includegraphics[angle=0, width=0.9\textwidth]{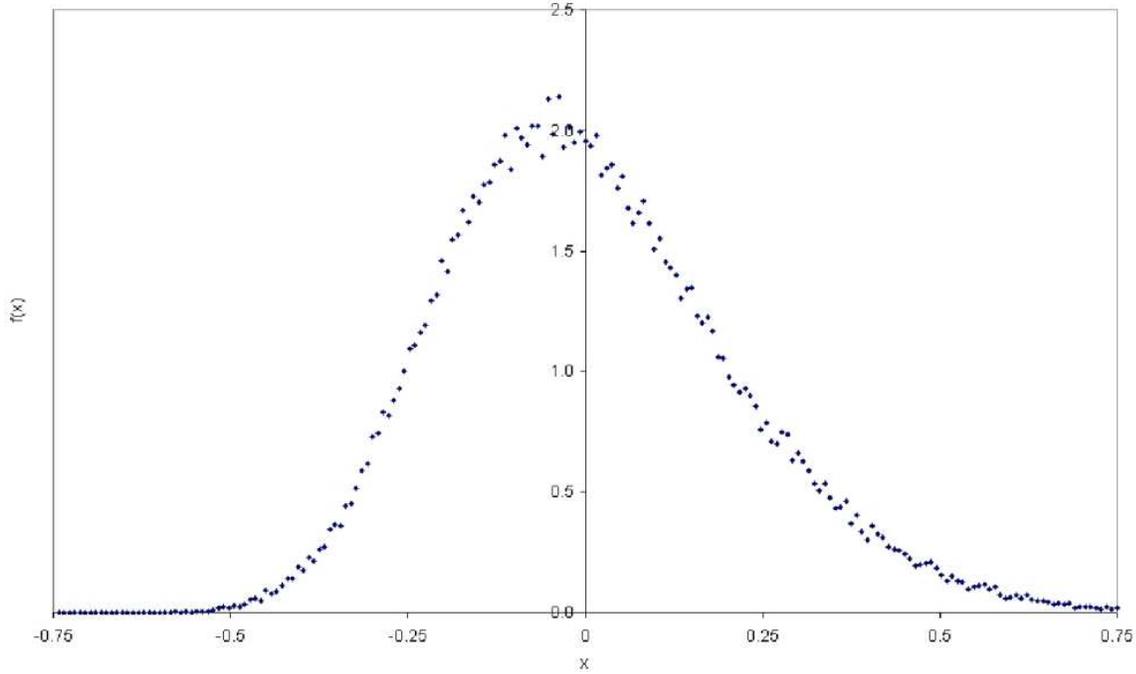}
\end{center}
\caption{Estimated probability density function for $\tGo$.}
\label{pdffig}
\end{figure}

(e) Figure \ref{pdffig} is a plot of the estimated probability
density function of $\tGo$ given by
(\ref{1006b}). This was obtained by
performing $10^5$ repeated
simulations of the ONG on a sequence of
$10^3$ uniform (simulated) random points on $(0,1)$. For each
simulation, the expected value
of $\OO^{1,1}(\U_{10^3})$
was subtracted
from the total length of the simulated ONG
to give an approximate realization of the
distributional limit.
The density function was then estimated from
the sample of $10^5$ realizations.
The simulated
sample from which the density
estimate
was taken had sample mean $\approx 3 \times 10^{-3}$
and sample variance $\approx 0.0425$, which are reasonably
close to the expectation and variance of $\tGo$.

\subsection{The nearest-neighbour (directed) graph}
\label{secnng}

Our next result gives exact expressions for the expectation and
variance of the total weight of the the
nearest-neighbour (directed) graph on $n$ independent uniform
random points in the unit interval. The
nearest-neighbour (directed) graph on a point set $\X$ places a directed edge from each vertex
to its nearest-neighbour (in the Euclidean sense). 

Let $\LL_1^{1,\alpha}(\X)$ denote the
total weight, with weight function $w_\alpha$ given by (\ref{wf}), of the nearest-neighbour (directed) graph
on vertex set $\X \subset (0,1)$. We use this notation to be consistent with \cite{llnpaper}, which
presents explicit laws of large numbers for nearest-neighbour graphs including this one. Let $\U_n$ denote
the binomial point process consisting of 
$n$ independent uniform random points in the unit interval. In this section with
give explicit results for the expectation and variance
of $\LL_1^{1,\alpha}(\U_n)$.

Let$~_2
F_1(\cdot,\cdot;\cdot;\cdot)$ denote the Gauss hypergeometric
function (see e.g.~Chapter 15 of \cite{as}) defined for $|z|<1$ and $c\neq 0,-1,-2,\ldots$ by \bea
\label{gauss} _2F_1(a,b;c;z) := \sum_{i=0}^\infty
\frac{(a)_i(b)_i}{(c)_n i!} z^i ,\eea where $(a)_i$ is
Pochhammer's symbol $(a)_i:=\Gamma(a+i)/\Gamma(a)$.
 For $n \in \{ 2,3,\ldots\}$, $\alpha>0$, set
\bea \label{0530c} J_{n,\alpha} :=
6^{-\alpha-1} \frac{\Gamma(n+1) \Gamma(2+2\alpha)}{(1+\alpha)\Gamma(n+1+2\alpha)}
~_2 F_1(-\alpha,1+\alpha;2+\alpha;1/3).
\eea
Also, for $\alpha>0$, set
\bea
\label{0530e}
j_\alpha := 8\lim_{n \to \infty} ( n^{2\alpha} J_{n,\alpha} )=
8 \cdot 6^{-\alpha-1} \frac{\Gamma(2+2\alpha)}{1+\alpha}~_2 F_1(-\alpha,1+\alpha;2+\alpha;1/3).
\eea
\begin{theorem}
\label{nng1d}
Suppose $\alpha>0$. For $n \in \{ 2,3,4,\ldots\}$ we have
\bea
\label{0530a}
\Exp [ \LL_1^{1,\alpha} (\U_n)] = ((n-2)2^{-\alpha}+2) \frac{ \Gamma(n+1)
\Gamma(\alpha+1)}{\Gamma(n+\alpha+1)} \sim 2^{-\alpha} \Gamma(\alpha+1) n^{1-\alpha},
\eea
as $n \to \infty$. Also, for $n \in \{ 4,5,6,\ldots\}$ 
\bea
\label{0530b}
\Var[ \LL_1^{1,\alpha} (\U_n)] & = &
\frac{\Gamma(n+1)}{\Gamma (n+2\alpha+1)}
\left[ \Gamma(2\alpha+1) (2-2 \cdot 3^{-2\alpha}+4^{-\alpha}n+2 \cdot 3^{-1-2\alpha}n) \right. \nonumber\\
&  & \left. + \Gamma (\alpha+1)^2 ( 4 + 12 \cdot 4^{-\alpha} - 12 \cdot 2^{-\alpha}
+ 2^{2-\alpha} n - 7 \cdot 4^{-\alpha} n +4^{-\alpha} n^2 ) \right] \nonumber\\
&  & - \left( \Exp [ \LL_1^{1,\alpha} (\U_n)] \right)^2 + 8(n-3) J_{n,\alpha},
\eea
where $\Exp [ \LL_1^{1,\alpha} (\U_n)]$
is given by (\ref{0530a}) and $J_{n,\alpha}$ is given by (\ref{0530c}).
Further, for $\alpha>0$
\bea
\label{0530d}
n^{2\alpha -1} \Var[ \LL_1^{1,\alpha} (\U_n)] \to
(4^{-\alpha} +2 \cdot 3^{-1-2\alpha} ) \Gamma (2\alpha+1) - 4^{-\alpha}
(3+\alpha^2) \Gamma (\alpha+1)^2  + j_\alpha ,\eea
as $n \to \infty$, where $j_\alpha$ is given by (\ref{0530e}).
\end{theorem}
Using (\ref{0530b}), with (\ref{0530c}), one obtains, for instance
\[ \Var [ \LL_1^{1,1} (\U_n) ] = \frac{2n^2 +17n+12}{12(n+1)^2 (n+2)} = \frac{1}{6} n^{-1} + O(n^{-2}),\]
and
\[ \Var [ \LL_1^{1,2} (\U_n) ] = \frac{85n^3 +3645n^2+7154n-456}{108(n+1)^2 (n+2)^2 (n+3)(n+4)}
 = \frac{85}{108} n^{-3} + O(n^{-4}).\]
Also,
the limiting constants $j_\alpha$ can be evaluated explicitly, so that one
can obtain values for $V_\alpha :=
\lim_{n \to \infty} ( n^{2\alpha -1} \Var[ \LL_1^{1,\alpha} (\U_n)])$.
Table \ref{tabvar} below gives some values of $V_\alpha$. 
We prove Theorem \ref{nng1d} in Section \ref{nngvar}.
\begin{table}[!h]
\begin{center}
\begin{tabular}{|l|l|l|l|l|l|}
\hline
$\alpha$  &  $\frac{1}{2}$ & 1 & 2 & 3 & 4   \\
\hline
$V_\alpha$ &   $\frac{1}{2} +\sqrt{2} \arcsin \left( \frac{1}{\sqrt{3}} \right)
-\frac{13\pi}{32} \approx 0.094148$ & $\frac{1}{6}$ &  $\frac{85}{108}$ & $\frac{149}{18}$ & $\frac{135793}{972}$  \\
\hline
\end{tabular}
\caption{Some values of $V_\alpha$.}
\label{tabvar}
\end{center}
\end{table}
One can obtain analogous explicit results in the case of $\LL_1^{1,\alpha} (\Po_n)$, where
$\Po_n$ is a homogeneous Poisson point process of intensity $n$
 on $(0,1)$: see \cite{thesis}, where a ``Poissonized''
version of (\ref{0530d}) is given.

The remainder of the present paper is organized as follows. Our results on the ONG in general dimensions
(Theorem \ref{onngthm} (ii) and (iii)) are proved in Section \ref{prf1027}. The main body of this
paper, Section \ref{1d}, is devoted to the ONG in one dimension and the proof of Theorem \ref{onng1}. In Section \ref{nngvar} we prove Theorem \ref{nng1d}. Finally, in the Appendix, we give the proofs of
some technical lemmas which would
otherwise interrupt the flow of the paper.

\section{Proof of Theorem \ref{onngthm} (ii) and (iii)}
\label{prf1027}

Suppose $d\in \N$. For $i\in\N$, let $Z_i(d) := \OO^{d,1} (\U_i) - \OO^{d,1} (\U_{i-1})$,
setting $\OO^{d,1} (\U_{0}) :=0$. That is, $Z_i(d)$ is the gain in length
of the $\onng$ on a sequence of independent uniform random points in $(0,1)^d$
on the addition of the $i$th point. Let $d_1(\bx;\X)$ denote the (Euclidean)
distance between $\bx \in \R^d$ and its nearest-neighbour in the
point set $\X \subset \R^d$.

\begin{lemma}
\label{1027e}
For $\alpha>0$ and $d \in \N$, as $n \to \infty$,
\bea
\label{1027f}
\Exp [ (Z_n(d)) ^\alpha] = O (n^{-\alpha/d}).
\eea
\end{lemma}
\proof  
We have
\[ \Exp [ (Z_n(d))^\alpha]
= \Exp[ (d_1 (\bU_1 ; \U_n))^\alpha]=
n^{-\alpha/d} \Exp[ (d_1 (n^{1/d} \bU_1 ; n^{1/d} \U_n ))^\alpha],\]
which is $O(n^{-\alpha/d})$ (see the proof of Lemma 3.3 in \cite{llnpaper}). $\square$ \\

\rem
We can obtain, by some more detailed
analysis, (see \cite{thesis})
\bean
\Exp [ (Z_n(d)) ^\alpha] = \frac{\alpha}{d} ( nv_d)^{-\alpha/d} \Gamma (\alpha/d) + o(n^{-(\alpha/d)}) .
\eean

\noindent
{\bf Proof of Theorem \ref{onngthm} (ii) and (iii).}
With the definition of $Z_i(d)$ in this section, let
\bean
W(d,\alpha) = \sum_{i=1}^\infty (Z_i(d))^\alpha.
\eean
The sum converges almost surely since it has non-negative terms
and, by (\ref{1027f}), has finite expectation for $\alpha>d$.
Let $k \in \N$. By (\ref{1027f}) and H\"older's inequality, there exists a constant $C\in(0,\infty)$ such that
\bean \Exp [ (W(d,\alpha))^k ] & = & \sum_{i_1=1}^\infty \sum_{i_2=1}^\infty \cdots
\sum_{i_k=1}^\infty  \Exp[ (Z_{i_1}(d))^\alpha (Z_{i_2}(d))^\alpha \cdots (Z_{i_k}(d))^\alpha]
\\ & \leq & C \sum_{i_1=1}^\infty \sum_{i_2=1}^\infty \cdots
\sum_{i_k=1}^\infty  i_1^{-\alpha/d} i_2^{-\alpha/d} \cdots i_k^{-\alpha/d} < \infty,\eean
since $\alpha/d >1$. The $L^k$ convergence then follows from the dominated convergence
theorem, and we have part (ii) of Theorem \ref{onngthm}.

Finally, for (iii) of Theorem \ref{onngthm}, we have, when $\alpha=d$
\[ n (Z_n(d))^d \eqd (d_1 (n^{1/d} \bU_1 ; n^{1/d} \U_n))^d \tod d_1(\0; \H_1)^d,\]
by the proof of Lemma 3.2 of \cite{penyuk2}. Since the sequence $(d_1 (n^{1/d} \bU_1 ; n^{1/d} \U_n))^d$
is uniformly integrable (see the proof of Theorem 2.4 of \cite{penyuk2}) we have
\[ \Exp[ n (Z_n(d))^d ] \to \Exp [ (d_1( \0 ; \H_1))^d] = v_d^{-1} ,\]
where the last inequality follows by a simple computation, or
by equation (2.7) of \cite{llnpaper}. So $\Exp[ (Z_n(d))^d ] =n^{-1}(v_d^{-1}+h(n))$ where $h(n) \to 0$
as $n \to \infty$. Thus
\[ \Exp \sum_{i=1}^n (Z_i(d))^d = \sum_{i=1}^n i^{-1}(v_d^{-1}+h(i)) = v_d^{-1} \log{n} +o(\log{n}),\]
and so we have (\ref{1029a}),
completing the proof of Theorem \ref{onngthm}. $\square$

\section{The $\onng$ in $d=1$} \label{1d}

\subsection{Notation and results}
\label{notres}
In this section we analyse
the $\onng$
in the interval $(0,1)$.
Theorem \ref{onng1} will follow from the main result of this section,
Theorem \ref{onng01} below. We introduce our notation.

For any finite sequence of points $\TT_n = ( x_1,x_2,\ldots,x_n ) \in
[0,1]^n$ with distinct inter-point distances, we construct the $\onng$ as follows.
Insert the points $x_1,x_2,\ldots$ into $[0,1]$ in order, one at a
time.
We join
a new point by an edge to its nearest neighbour among those already
present, provided that such a point exists. In other words, for
each point $x_i$, $i \geq 2$, we join $x_i$
by an edge to the point of $\{ x_j : 1 \leq j < i \}$ that
minimizes $|x_i-x_j|$. In this way we
construct a tree rooted at $x_1$,
which we denote by $\onng
( \TT_n )$. Denote the total weight (under
weight function $w_\alpha$ given by (\ref{wf}),
 $\alpha>0$) of $\onng (\TT_n)$ by
 $\OO^{1,\alpha}
( \TT_n )$, to be consistent with our previous
notation.

For what follows, our main interest is the
case in which $\TT_n$ is a random vector in $[0,1]^n$.
In this case,
set $\tO^{1,\alpha} (
\TT_n ) := \OO^{1,\alpha} ( \TT_n ) -
\Exp [ \OO^{1,\alpha} ( \TT_n ) ]$, the
centred total weight of the $\onng$ on $\TT_n$.
Let $(U_1,U_2,U_3,\ldots)$ be a sequence
of independent uniformly distributed random
variables in $(0,1)$, and for $n \in \N$
set $\U_n := (U_1,U_2, \ldots, U_n)$.
Given $\U_n$, we define the augmented sequences
$\U_n^0 = (0,U_1,\ldots,U_n)$
and $\U_n^{0,1} = (0,1,U_1,\ldots,U_n)$.
Notice that
 $\onng ( \U_n^{0,1} )$
and $\onng (\U_n^0)$ both give a tree rooted at $0$,
and that in $\onng (\U_n^{0,1})$ the first edge
is from $1$ to $0$.

We now state the main result of this section, from which Theorem
\ref{onng1} will follow. The convergence of joint distribution results in (\ref{1703b}) and
(\ref{0202c}) are given in more detail, complete with joint distribution fixed-point
representation, in Propositions \ref{1003e} and \ref{0201a}.

\begin{theorem}
\label{onng01}
\begin{itemize}
\item[(i)]
For $1/2<\alpha<1$, we have that, as $n \to \infty$,
\bea
\label{1703b}
(\tO^{1,\alpha} (\U_n^{0,1}),
\tO^{1,\alpha} (\U_n^{0}),
\tO^{1,\alpha} (\U_n))
 \tod (\tJal,\tHal,\tGal),
\eea
where
$\tJal$, $\tHal$, $\tGal$ are jointly distributed random
variables with marginal distributions given by (\ref{0919z}),
(\ref{0923n}), (\ref{1006a}) respectively.
\item[(ii)]
For $\alpha =1$, we have that, as $n \to \infty$,
\bea
\label{0202c}
( \tO^{1,1} (\U_n^{0,1}),
\tO^{1,1} (\U_n^{0}),
\tO^{1,1} (\U_n))
\tod (\tJo,\tHo,\tGo),
\eea
where
$\tJo$, $\tHo$, $\tGo$  are jointly distributed random
variables with marginal distributions given by (\ref{0919x}),
(\ref{0924e}), (\ref{1006b}) respectively.
The first three moments of $\tJo$, $\tHo$ and $\tGo$ are given in Table
\ref{tabmoms}. Further, the variables on the right hand side of (\ref{0202c})
satisfy $\Cov(\tJo,\tHo)= ((9+6\log{2})/32)-(\pi^2/24) \approx -1.84204 \times 10^{-5}$,
$\Cov(\tGo,\tHo)= ((35+10\log{2})/48)-(\pi^2/24) \approx 0.0255536$, and
$\Cov(\tGo,\tJo)=((7+4\log{2})/24)-(\pi^2/24) \approx -4.04232 \times 10^{-3}$.
\item[(iii)]
For $\alpha > 1$, we have that,
as $n \to \infty$,
\[ \OO^{1,\alpha} (\U_n^{0,1}) \to 1+\Jal; ~~~\OO^{1,\alpha} (\U_n^{0}) \to \Hal ,\]
where the convergence is almost sure
and in $L^p$, $p\in \N$, and the distributions
of $\Jal$ and $\Hal$ are
 given by (\ref{0919y}) and (\ref{0923l}) respectively.
\end{itemize}
\end{theorem}

\begin{table}[!h]
\begin{center}
\begin{tabular}{|l|l|l|l|}
\hline
 &  $\Exp[\cdot]$ & $\Var[\cdot]$ & $\Exp[(\cdot)^3]$   \\
\hline
$\tJo$  &  0 & $((1+\log{2})/4)-(\pi^2/24) \approx 0.012053$ & $\approx -0.00005733$ \\
$\tHo$  &  0 & $((3+\log{2})/8)-(\pi^2/24) \approx 0.050410$ & $\approx 0.00323456$\\
$\tGo$  &  0 & $((19+4\log{2})/48)-(\pi^2/24) \approx 0.042362$ &
$\approx 0.00444287$ \\
\hline
\end{tabular}
\caption{First three moments for the random variables $\tJo$, $\tHo$, $\tGo$.}
\label{tabmoms}
\end{center}
\end{table}

Our method for establishing convergence in distribution results
is based on the {\em recursive} nature of the $\onng$. Essential
is its {\em self-similarity} (scaling
property). In terms of the total weight,
this says that for any $t \in (0,1)$, if $V_1,\ldots,V_n$
are independent and uniformly distributed on $(0,t)$,
then the distribution of $\OO^{1,\alpha} (V_1, \ldots, V_n)$
is the same as that of $t^\alpha \OO^{1,\alpha} ( U_1, \ldots, U_n)$.

Write $U=U_1$ for the position of the first arrival.
For ease of notation, denote \bea
\label{0999}
Y_n := \OO^{1,\alpha} (\U_n^{0,1})-1,
\eea
where by subtracting 1 we discount the length of the edge from 1 to 0.
Then
using
the self-similarity of the $\onng$, and conditioning
on the first arrival,
we have the following relations:
\bea
\label{0926ii}
\OO^{1,\alpha} (\U_n) \eqd U^\alpha \OO^{1,\alpha}_{\{1\}} (\U^0_{N(n)})
+ (1-U)^\alpha \OO^{1,\alpha}_{\{2\}} (\U^0_{n-1-N(n)}), \\
 \OO^{1,\alpha}(\U_n^0) \eqd
U^\alpha \OO^{1,\alpha}_{\{1\}} (\U_{N(n)}^{0,1})
 +(1-U)^\alpha \OO^{1,\alpha}_{\{2\}} (\U_{n-1-N(n)}^0),
\label{0923m} \\
\label{1029f}
 Y_n \eqd (\min \{U,1-U\})^\alpha + U^\alpha Y^{\{1\}}_{N(n)}
+ (1-U)^\alpha Y^{\{2\}}_{n-1-N(n)} ,
\eea
where, given $U$, $N(n) \sim {\rm Bin} (n-1,U)$ gives the number
of points of $U_2,U_3,\ldots,U_n$ that arrive to the left
of $U_1=U$. Given $U$ and $N(n)$,
 $\OO^{1,\alpha}_{\{1\}} (\cdot)$ and
$\OO^{1,\alpha}_{\{2\}} (\cdot)$ are
independent copies of $\OO^{1,\alpha} (\cdot)$. Also, given $U$ and $N(n)$, $Y_{N(n)}^{\{1\}}$ and
 $Y_{n-1-N(n)}^{\{2\}}$
are independent with the distribution of $Y_{N(n)}$ and
 $Y_{n-1-N(n)}$, respectively.

For $\alpha>1$, we prove almost sure and $L^p$ $(p \in \N)$ convergence of $\OO^{1,\alpha} (\U_n^0)$
and $\OO^{1,\alpha} (\U_n^{0,1})$, in the same way as in the proof
of Theorem \ref{onngthm} (ii), and thereby obtain the corresponding result for
$\OO^{1,\alpha} (\U_n)$. The relations (\ref{0926ii}), (\ref{0923m}) and (\ref{1029f})
will then enable us to prove the desired results for $\alpha>1$.

For $1/2<\alpha \leq 1$, we use a result of Neininger and
R\"uschendorf \cite{neinrusch} on limit theorems for `divide and
conquer' recurrences. However, we cannot apply this directly to
(\ref{0926ii}) to obtain the convergence of $\OO^{1,\alpha}
(\U_n)$, since (\ref{0926ii}) is not of the required form; the
variables on the right are not of the same type as the variable on
the left. On the other hand, we see that (\ref{1029f}) is of the
desired form. This will be the basis of our analysis for $1/2 <
\alpha \leq 1$.

Indeed, by considering a vector defined in terms of all three of
$\OO^{1,\alpha} (\U_n)$, $\OO^{1,\alpha} (\U_n^0)$, and $\OO^{1,\alpha} (\U_n^{0,1})$,
we obtain the recurrence relation (\ref{0920iq}) below.
We can then apply the result of \cite{neinrusch}.
This is why we need to consider $\OO^{1,\alpha} (\U_n^0)$ and $\OO^{1,\alpha} (\U_n^{0,1})$
in addition to $\OO^{1,\alpha} (\U_n)$.

The outline of the remainder of this section is as follows. In
Section \ref{spacings} below, we give a discussion of the theory
of spacings, which will be very useful in the sequel. In Section
\ref{sec0} we begin our analysis of the $\onng$ with some
preliminary results, based on the discussion in Section
\ref{spacings}. Then, in Sections  \ref{sec2}, \ref{sec1} and
\ref{sec3} we give results on $\OO^{1,\alpha} (\cdot)$ when
$1/2<\alpha<1$, $\alpha=1$, and $\alpha>1$ respectively. Finally, in
Section \ref{sec4} we give a proof of Theorems \ref{onng01} and
\ref{onng1}.

\subsection{Spacings}
\label{spacings}

The one-dimensional
models considered in this paper (the $\onng$ and the standard
nearest-neighbour graph) are defined in terms of the {\em spacings}
of points in the unit interval. Thus the
theory of so-called Dirichlet spacings will be
useful. For some
general references on spacings, see for example \cite{pyke}. A large number of statistical
tests are based on spacings, see e.g.~\cite{darling} for a few examples.

Recall that $\U_n$ denotes the binomial point process
consisting of $n$ independent
uniform random variables on $(0,1)$, $U_1, U_2, \ldots, U_n$.
Given $\{ U_1,\ldots, U_n \} \subseteq (0,1)$, denote
the order statistics of $U_1,\ldots,U_n$, taken in increasing order,
as $U^n_{(1)},U_{(2)}^n,\ldots,U_{(n)}^n$. Thus $(U^n_{(1)},\ldots,U_{(n)}^n)$
is a nondecreasing sequence, forming a permutation of the original $(U_1,\ldots, U_n)$.

The points $U_1,\ldots,U_n$ divide $[0,1]$ into $n+1$ intervals. Denote the intervals
between points by $I_j^n:=(U^n_{(j-1)},U^n_{(j)})$ for $j=1,2,\ldots,n+1$, where we
set $U^n_{(0)}:=0$ and $U^n_{(n+1)}:=1$. Let the widths of these intervals (the spacings) be
\[ S^n_j := | I^n_j | = U^n_{(j)} - U^n_{(j-1)},\]
for $j=1,2,\ldots,n+1$.
For $n \in \N$, let $\Delta_n \subset \R^n$ denote the $n$-dimensional simplex,
that is
\bean
 \Delta_n:= \left\{ (x_1,\ldots,x_n) \in \R^n : x_i \geq 0, ~1 \leq i \leq n;~ \sum_{i=1}^{n} x_i \leq 1 \right\}.\eean
By the definition of $S_j^n$, we have that $S_j^n \geq 0$ for $j=1,\ldots,n+1$ and
$\sum_{j=1}^{n+1} S_j^n =1$. So we see that the vector $(S_1^n, S_2^n, \ldots, S_{n+1}^n)$
is completely specified by any $n$ of its $n+1$ components, and any such $n$-vector
belongs to the simplex $\Delta_n$. It is not hard to show that any such $n$-vector is,
in fact, uniformly distributed over the
simplex.
Hence $( S_1^n,
\ldots , S_n^n )$ is uniform over the simplex $\Delta_n$, and
$S_{n+1}^n = 1 -\sum_{i=1}^n S_i^n$.

Thus $(S_1^n, S_2^n, \ldots, S_{n+1}^n)$ has the
symmetric {\em Dirichlet distribution} with parameter $1$ (see, e.g.,~\cite{bill}, p.~246),
and any $n$-vector of the $S_j^n$ has the Dirichlet density
\bea
\label{0606a}
 f(x_1,\ldots,x_n) = n!, ~~~ (x_1,\ldots,x_n) \in \Delta_n .\eea
In particular, the spacings $S_j^n$, $j=1,\ldots,n+1$ are {\em exchangeable} -- the distribution
of  $(S_1^n, S_2^n, \ldots, S_{n+1}^n)$ is invariant under any permutation of its components.

By integrating out over the simplex, from (\ref{0606a}) one can readily obtain the marginal
distributions for the spacings. Thus, for $n \geq 1$,
a single spacing has density
\bea
\label{0606b}
f(x_1) = n(1-x_1)^{n-1}, ~~~ 0 \leq x_1 \leq 1,
\eea
while for $n \geq 2$, any two spacings have joint
density
\bea
\label{0606c}
f(x_1,x_2) = n(n-1)(1-x_1-x_2)^{n-2}, ~~~ (x_1,x_2) \in \Delta_2,
\eea
and for $n \geq 3$ any three spacings have joint density
\bea
\label{0606d}
f(x_1,x_2,x_3) = n(n-1)(n-2)(1-x_1-x_2-x_3)^{n-3}, ~~~ (x_1,x_2,x_3) \in \Delta_3.
\eea
Using the fact that (see, e.g., 6.2.1 in \cite{as})
\bea
\label{gint}
\int_0^1 t^{a-1} (1-t)^{b-1} \ud t = \frac{\Gamma(a)\Gamma(b)}{\Gamma(a+b)}
,\eea 
for $a>0$, $b>0$,
it then follows from (\ref{0606b}) that, for $\beta>0$, $n \geq 1$
\bea
\label{0530g}
\Exp \left[ (S_1^n)^\beta \right]
= \frac{\Gamma(n+1) \Gamma(\beta+1)}{\Gamma(n+\beta+1)},
\eea
and from (\ref{0606c}) that for $\beta>0$, $n \geq 2$
\bea
\label{0516b}
\Exp \left[ (S_1^n)^\beta (S_2^n)^\beta \right]
= \frac{\Gamma(n+1) \Gamma(\beta+1)^2}{\Gamma(n+2\beta+1)}.
\eea
When considering our nearest-neighbour graphs, 
we will encounter the {\em minimum} of two (or more) spacings. The following
results will also be needed in Section \ref{nngvar}.
\begin{lemma}
\label{lem0606}
For $n \geq 1$,
\bea
\label{1105b}
 \min \{ S_1^n, S_2^n \} \eqd S_1^n/2 .\eea
For $n \geq 2$,
\bea
\label{0516c}
 (S_1^n ,  \min \{S_2^n, S_3^n \} ) \eqd (S_1^n , S_2^n/2) .\eea
Finally, for $n \geq 3$
\bea
\label{0516d}
 ( \min \{ S_1^n, S_2^n\} , \min \{ S_3^n , S_4^n \} ) \eqd ( S_1^n/2 , S_2^n/2) ,\eea
and
\bea
\label{0516dd}
\min \{ S_1^n, S_2^n, S_3^n \} \eqd S_1^n/3.
\eea
\end{lemma}
\proof
We give the proof of (\ref{1105b}). The other results follow by very similar calculations
based on (\ref{0606c}) and (\ref{0606d}). Suppose $n \geq 2$.
From (\ref{0606c}), we have, for $0 \leq r \leq 1/2$
\bean
 \Pr [ \min \{S_1^n , S_2^n \} > r ] & = &
\Pr [ S_1^n >r,~ S_2^n>r] \\
& = & n(n-1) \int_r^{1-r} \ud x_1 \int_r^{1-x_1} (1-x_1-x_2)^{n-2} \ud x_2 \\
& = & (1-2r)^n = \Pr [ S_1^n > 2r ],\eean
and so we have (\ref{1105b}). $\square$

\subsection{Preparatory results}
\label{sec0}

We now return to the $\onng$.
We make use of the discussion of spacings in Section \ref{spacings}.
For $n\in \N$
let $Z_n$, $H_n$ and $T_n$ denote the
random variables given by the gain in length, on the addition of the point $U_n$, of the $\onng$
on $\U_{n-1}$, $\U_{n-1}^0$ and $\U_{n-1}^{0,1}$ respectively.
That is, with the convention $\OO^{1,1} (\U_0) = \OO^{1,1} (\U_0^0)=0$
and $\OO^{1,1}(\U_0^{0,1})=1$, for $n\in\N$ set
\bea
\label{0918b}
Z_n :=
\OO^{1,1}(\U_n)-\OO^{1,1}(\U_{n-1}), \\
H_n :=
\OO^{1,1}(\U_n^{0})-\OO^{1,1}(\U_{n-1}^{0}), \nonumber\\
T_n :=
\OO^{1,1}(\U_n^{0,1})-\OO^{1,1}(\U_{n-1}^{0,1}). \nonumber
\eea
Thus, for example, in the $\onng(\U_n^{0,1})$
with weight function $w_\alpha$ as given by (\ref{wf}),
the $n$th edge
to be added has weight $T_n^\alpha$.

We will make use of the following discussion for the proof of Lemma \ref{0204a} below.
For $\alpha>0$, with the definitions at (\ref{0918b}), we have that
\bea
\label{0204d}
\OO^{1,\alpha} (\U_n^{0}) - \OO^{1,\alpha} (\U_n^{0,1})
& = & -1+\sum_{i=1}^n \left(H_i^\alpha - T_i^\alpha \right), ~~{\rm and} \\
\OO^{1,\alpha} (\U_n) - \OO^{1,\alpha} (\U_n^0)
& = &
\sum_{i=1}^n \left(Z_i^\alpha - H_i^\alpha \right)
= -H_1^\alpha + \sum_{i=2}^n \left(Z_i^\alpha - H_i^\alpha \right),
\eea
since $Z_1=0$.
Consider the arrival of the point $U_n$. For any $n$, $T_n$ and $H_n$ are the same unless the point $U_n$
falls in the right hand half of the rightmost interval $I^{n-1}_n$ of width $S^{n-1}_n$. Denote this latter event
by $E_n$. Given $S^{n-1}_n$,
the probability of $E_n$ is $S^{n-1}_n/2$. Given $S^{n-1}_n$, and given that $E_n$ occurs,
the value of $T_n$ is given by $(1-V_n) S^{n-1}_n/2$ and the value of $H_n$ by
$(1+V_n) S^{n-1}_n/2$, where $V_n=1+2(U_n -1)/S^{n-1}_n$ is uniform on $(0,1)$ given $E_n$. So we have
that, for $n \in \N$,
given $S_n^{n-1}$
\bea
\label{0204e}
 H_n^\alpha - T_n^\alpha =  \1_{E_n} \left( \frac{S^{n-1}_n}{2} \right)^\alpha \left( (1+V_n)^\alpha
-(1-V_n)^\alpha \right),\eea
where $E_n$ is an event with probability $S^{n-1}_n/2$.
A similar argument (based this time on the {\em
leftmost} spacing) yields that, for $n \geq 2$
\bea
\label{0204ff}
 Z_n^\alpha - H_n^\alpha =  \1_{F_n} \left( \frac{S^{n-1}_1}{2} \right)^\alpha \left( (1+W_n)^\alpha
-(1-W_n)^\alpha \right),\eea
where $F_n$ is an event with probability $S^{n-1}_1/2$ and, given $F_n$, $W_n$ is uniform
on $(0,1)$. 

We will need the following asymptotic expansion, which follows from Stirling's formula (see e.g. 6.1.37 in \cite{as}). For any $\beta>0$, as $n \to \infty$,
\bea
\label{stir}
\frac{\Gamma(n+1)}{\Gamma(n+1+\beta)} = n^{-\beta} - \frac{1}{2} \beta (\beta+1) n^{-\beta-1} + O(n^{-\beta-2}).
\eea
\begin{lemma}
\label{0204a}
For $\alpha>0$ and $n \geq 2$,
we have that
\bea
\label{0204c}
\Exp[\OO^{1,\alpha} (\U_n^{0}) - \OO^{1,\alpha} (\U_n^{0,1})]
& = & \frac{1-2^{-\alpha}-\alpha}{\alpha} + (2^{-\alpha}-1)\frac{\Gamma(\alpha)\Gamma(n+1)}{\Gamma(n+1+\alpha)} \nonumber\\
& = & \frac{1-2^{-\alpha}-\alpha}{\alpha} +O(n^{-\alpha}) ,\eea
and
\bea
\label{0204cc}
\Exp[\OO^{1,\alpha} (\U_n) - \OO^{1,\alpha} (\U_n^0)]
& = & \frac{1-2^{-\alpha}-\alpha}{\alpha(1+\alpha)} +(2^{-\alpha}-1) \frac{\Gamma(\alpha)\Gamma(n+1)}{\Gamma(n+1+\alpha)}  \nonumber\\
& = & \frac{1-2^{-\alpha}-\alpha}{\alpha(1+\alpha)} +O(n^{-\alpha}).\eea
\end{lemma}
\proof Suppose $\alpha>0$.
From
(\ref{0204e}) we have that for $n \in\N$
\[ \Exp[H_n^\alpha-T_n^\alpha | S^{n-1}_n] =
 (S^{n-1}_n)^{1+\alpha} \left( \frac{1-2^{-\alpha}}{1+\alpha}\right).\]
 So by (\ref{0530g}) we have that
\bean
\Exp[H_n^\alpha-T_n^\alpha ] =
 \frac{(1-2^{-\alpha}) \Gamma(1+\alpha) \Gamma(n)}{\Gamma(n+1+\alpha)}.\eean 
Thus, from (\ref{0204d}),
\[ \Exp[\OO^{1,\alpha} (\U_n^{0}) - \OO^{1,1} (\U_n^{0,1})]
= -1+ \Exp \sum_{i=1}^n (H^\alpha_i-T^\alpha_i) =
-1+ \sum_{i=1}^n \frac{(1-2^{-\alpha}) \Gamma(1+\alpha) \Gamma(i)}{\Gamma(i+1+\alpha)},\]
the last equality following by induction on $n$. This then gives (\ref{0204c}), 
with the asymptotic expression following by (\ref{stir}).
Similarly, from (\ref{0204ff})
\bean
\Exp[Z_n^\alpha-H_n^\alpha ] =
 \frac{(1-2^{-\alpha}) \Gamma(1+\alpha) \Gamma(n)}{\Gamma(n+1+\alpha)},\eean
 for $n \geq 2$, while $\Exp[H_1^\alpha]=\Exp[U_1^\alpha]=(\alpha+1)^{-1}$ and
 $Z_1=0$.
With (\ref{stir}), (\ref{0204cc}) follows.  $\square$

\begin{lemma} \label{1215v}
\begin{itemize} \item[(i)]
For $n\in\N$, $T_n$ as defined at (\ref{0918b}) has distribution
function $F_n$ given by $F_n(t)=0$ for $t<0$, $F_n(t)=1$ for
$t \geq 1/2$, and
$F_n(t) = 1- ( 1-2t)^n$ for $0 \leq t \leq 1/2$.
\item[(ii)] For $\beta >0$,
 \bea
\label{0918c}
 \Exp[T_n^\beta] = 2^{-\beta} \frac{
\Gamma(n+1) \Gamma(\beta +1 )}{\Gamma ( n+ \beta +1 )} . \eea
In particular,
\bea
\label{0918d}
\Exp[T_n] = \frac{1}{2(n+1)} ; ~~~ \Var[T_n]
= \frac{n}{4(n+1)^2(n+2)}.
\eea
\item[(iii)] For $\beta>0$, as $n \to \infty$
\bea
\label{0920d}
 \Exp[T_n^\beta] = 2^{-\beta} \Gamma( \beta +1) n^{-\beta} +O(n^{-\beta-1}).\eea
\item[(iv)] As $n \to \infty$,
\[ 2nT_n \tod
\mathrm{Exp}(1), \] where $\mathrm{Exp}(1)$ is an exponential
random variable with parameter 1. \end{itemize}
\end{lemma}
\proof
By conditioning on the number of $U_j$, $j \leq n$ with $U_j \leq U_n$, using
Lemma \ref{lem0606}, and by exchangeability of the spacings, we have that for
$n \geq 1$,
$T_n \eqd \min \{ S_1^n, S_2^n \} \eqd S_1^n/2$, by (\ref{1105b}). Then (i) follows
by (\ref{0606b}), and (ii) follows by (\ref{0530g}). 
Part (iii) then follows from part (ii)
by (\ref{stir}).
For (iv), we have that, for $t \in [0,\infty)$,
and $n$ large enough so that $t/(2n) \leq 1/2$,
 \bean
 \Pr [2nT_n > t] =
\Pr [ T_n > t/(2n) ] = \left( 1 - (t/n)
\right)^n  \to  e^{-t} ,
\eean
as $n \to \infty$,
but $1-e^{-t}$, $t \geq 0$
is the distribution function of an exponential random variable
with parameter 1.
$\square$

\begin{proposition}
\label{propw1}
Recall that $\gamma \approx 0.57721566$ is Euler's constant, defined at (\ref{gamma}).
Suppose $\alpha>0$. As $n \to \infty$, we have
\bea
\label{0919wa}
 \Exp[ \OO^{1,\alpha}(\U_n^{0,1}) ] & = &
\frac{\Gamma(\alpha+1)}{1-\alpha} 2^{-\alpha}n^{1-\alpha} + 1-\frac{2^{-\alpha}}{1-\alpha}
+O(n^{-\alpha}); ~~~  (0<\alpha<1) \\
\label{0919wb}
\Exp[ \OO^{1,1}(\U_n^{0,1}) ]  & = & \frac{1}{2}
\log{n} + \frac{1}{2}(\gamma+1) +O(n^{-1}); ~~~  \\
\label{0919wc}
\Exp[ \OO^{1,\alpha}(\U_n^{0,1}) ] & = & 1+\frac{2^{-\alpha}}{\alpha-1}
+O(n^{1-\alpha}) ~~~ (\alpha>1)
 \eea
\end{proposition}
\proof
Counting the first edge from $1$ to $0$,
we have \[ \Exp[\OO^{1,\alpha}(\U_n^{0,1})] = 1+\sum_{i=1}^n
\left( \Exp [ \OO^{1,\alpha}(\U_i^{0,1}) ] - \Exp[\OO^{1,\alpha}(\U_{i-1}^{0,1})]
\right) = 1+ \sum_{i=1}^n \Exp[T_i^\alpha] .
\] In the case where $\alpha=1$, $\Exp[T_i] = (2(i+1))^{-1}$ by
(\ref{0918d}), and (\ref{0919wb}) follows by (\ref{gamma}).
For general $\alpha>0$, $\alpha \neq 1$,
from (\ref{0918c}) we have that
\bea
 \label{1117a}
\Exp[\OO^{1,\alpha}(\U_n^{0,1})] & = & 1+2^{-\alpha} \Gamma(1+\alpha) \sum_{i=1}^n
\frac{\Gamma(i+1)}{\Gamma(1+\alpha+i)} \nonumber\\
& = & 1 + \frac{2^{-\alpha}}{\alpha-1} -
\frac{2^{-\alpha} \Gamma(1+\alpha) \Gamma(n+2)}{(\alpha-1) \Gamma(n+1+\alpha)}
,
 \eea 
 the final equality proved by induction on $n$.
 By Stirling's formula, the last term satisfies
\begin{equation} \label{1117b} -\frac{2^{-\alpha}
\Gamma(1+\alpha) \Gamma(n+2)}{(\alpha-1)
\Gamma(n+1+\alpha)} = -2^{-\alpha} \frac{\Gamma(1+\alpha)}{\alpha-1}
n^{1-\alpha} (1+O(n^{-1})), \end{equation} which tends to zero as
$n \to \infty$ for $\alpha >1$, to give us (\ref{0919wc}). For $\alpha <1$, we have (\ref{0919wa})
 from
(\ref{1117a}) and (\ref{1117b}).
$\square$

\begin{proposition}
\label{propw2}
Suppose $\alpha>0$. As $n \to \infty$, we have
\bea
\label{0923wa}
 \Exp[ \OO^{1,\alpha}(\U_n^0) ] & = &
\frac{\Gamma(\alpha+1)}{1-\alpha} 2^{-\alpha} n^{1-\alpha}  + \frac{1}{\alpha} - \frac{2^{-\alpha}}{\alpha(1-\alpha)}
+ O (n^{-\alpha}); ~~~  (0<\alpha<1) \\
\label{0923wb}
\Exp[ \OO^{1,1}(\U_n^0) ]   & = &
 \frac{1}{2}
\log{n} + \frac{1}{2}\gamma +O(n^{-1}); ~~~  \\
\label{0923wc}
\Exp[ \OO^{1,\alpha}(\U_n^0) ] & = & \frac{1}{\alpha}+\frac{2^{-\alpha}}{\alpha(\alpha-1)}
+O(n^{1-\alpha})~~~ (\alpha>1)
 \eea
\end{proposition}
\proof
This follows from Proposition \ref{propw1} with (\ref{0204c}).
$\square$

\subsection{Limit theory when $1/2<\alpha <1 $}
\label{sec2}

Let $U$ be uniform on $(0,1)$, and given $U$, let $N(n) \sim {\rm Bin} (n-1,U)$. Set
\bea
\label{0803d}
B_\alpha(n) :=
(n-1)^{1/2} \left( U^\alpha \left( \frac{N(n)}{n-1} \right)^{1-\alpha}
+ (1-U)^\alpha  \left( \frac{n-1-N(n)}{n-1} \right)^{1-\alpha} -1 \right)
.\eea

\begin{lemma}
\label{0803a}
Suppose $0 \leq \alpha \leq 1$.
Then, as $n \to \infty$,
\bea
\label{0903f}
B_\alpha(n)
\inLLL 0.
\eea
\end{lemma}
We defer the proof of this lemma to the Appendix. Note that for what follows in this paper we will only use $L^2$ convergence in (\ref{0903f}). However, the stronger $L^3$ version requires little extra work, and we will require the $L^3$ version
in future work dealing with the $\alpha \in (0,1/2]$ case.

\begin{proposition}
\label{1003e}
Suppose $1/2<\alpha<1$.
Then as $n \to \infty$,
\bea
\label{0217qq}
\left( \begin{array}{l} \tO^{1,\alpha} (\U_n^{0,1})
\\
\tO^{1,\alpha} (\U_n^{0}) - \tO^{1,\alpha} (\U_n^{0,1}) \\
\tO^{1,\alpha} (\U_n) - \tO^{1,\alpha} (\U_n^0)
\end{array}
\right)
\tod
\left( \begin{array}{l} \tJal
\\
\tR
\\
\tS
\end{array}
\right)
,\eea
where $(\tJal, \tR,\tS)$ satisfies the fixed-point equation
\bea
\label{0207aav}
\left( \begin{array}{l} \tJal \\ \tR  \\ \tS \end{array} \right)
& \eqd&
\left( \begin{array}{lll} U^\alpha & 0 & 0 \\ 0 & 0 & 0 \\ 0 & U^\alpha & 0 \end{array} \right)
\left( \begin{array}{l} \tJal^{\{1\}} \\ \tR^{\{1\}}
\\ \tS^{\{1\}} \end{array} \right)
+
\left( \begin{array}{lll} (1-U)^\alpha & 0 & 0 \\ 0 & (1-U)^\alpha & 0 \\ 0 & 0 & 0 \end{array} \right)
\left( \begin{array}{l} \tJal^{\{2\}} \\ \tR^{\{2\}} \\ \tS^{\{2\}}\end{array} \right) \nonumber\\ &  &
+
\left( \begin{array}{l} \min \{ U, 1-U \}^\alpha +\frac{2^{-\alpha}}{\alpha-1}((1-U)^\alpha +U^\alpha-1)
 \\ (U^\alpha-(1-U)^\alpha)  \1_{\{ U > 1/2 \}} +
 ((1-U)^\alpha-1) \frac{1-2^{-\alpha}}{\alpha} \\ (U^\alpha -\frac{1}{1+\alpha} )
 \frac{1-2^{-\alpha}-\alpha}{\alpha} \end{array}
 \right). \eea
In particular, $\tJal$ satisfies the fixed-point
equation (\ref{0919z}).
Also, $\Exp[ \tJal]=\Exp[\tR]=\Exp[\tS]=0$.
\end{proposition}
\proof
We make use
of Theorem 4.1 of \cite{neinrusch}, which is a general
result for `divide-and-conquer' type recurrences. Recall the definition
of $Y_n$ at (\ref{0999}). Let
\bea
\label{0995}
 R_n := \OO^{1,\alpha} (\U_n^{0}) - \OO^{1,\alpha} (\U_n^{0,1})+1, ~~~
S_n := \OO^{1,\alpha}(\U_n) - \OO^{1,\alpha} (\U_n^0). \eea
Write $U=U_1$ for the position
of the first arrival.
Given $U$, let $N(n) \sim {\rm Bin} (n-1,U)$ be the number
of points of $U_2,U_3,\ldots,U_n$ that arrive to the left
of $U_1=U$.
Using
the self-similarity of the $\onng$,
we have that $(Y_n,R_n,S_n)$ satisfies, for $\alpha>0$,
\bea
\label{0920iq}
\left( \begin{array}{l} Y_n \\ R_n
\\ S_n \end{array} \right)
 \eqd 
\left( \begin{array}{lll} U^\alpha & 0 & 0 \\ 0 & 0 & 0 \\ 0 & U^\alpha & 0 \end{array} \right)
\left( \begin{array}{l} Y^{\{1\}}_{N(n)} \\ R^{\{1\}}_{N(n)} \\
S^{\{1\}}_{N(n)} \end{array} \right) \nonumber\\
+
\left( \begin{array}{lll} (1-U)^\alpha & 0 & 0 \\ 0 & (1-U)^\alpha & 0 \\ 0 & 0 & 0 \end{array} \right)
\left( \begin{array}{l} Y^{\{2\}}_{n-1-N(n)} \\ R^{\{2\}}_{n-1-N(n)} \\
S^{\{2\}}_{n-1-N(n)} \end{array} \right)
 +
\left( \begin{array}{l} \min \{ U, 1-U \}^\alpha \\ (U^\alpha -(1-U)^\alpha) \1_{\{ U > 1/2 \}}
\\ -U^\alpha   \end{array} \right)
,
\eea
where, given $U$ and $N(n)$,
 $Y^{\{1\}}_{N(n)}$, $Y^{\{2\}}_{n-1-N(n)}$ are
independent copies of $Y_{N(n)}$, $Y_{n-1-N(n)}$ respectively, and similarly
for the $R$s and $S$s.
This equation
is of the form of (21) in \cite{neinrusch}. Suppose $1/2<\alpha<1$.
We now renormalise (\ref{0920iq}) by taking
\bea
\label{0998}
(\tY_n,\tR_n,\tS_n) := (Y_n - \Exp[Y_n],R_n-\Exp[R_n],S_n-\Exp[S_n]),
\eea so in the notation of \cite{neinrusch}, we
take $C_n \equiv 1$. That is,
\bea
\label{0996}
\tY_n = \tO^{1,\alpha} (\U_n^{0,1}),~~ \tR_n = \tO^{1,\alpha} (\U_n^{0}) - \tO^{1,\alpha} (\U_n^{0,1}),
~~
\tS_n = \tO^{1,\alpha} (\U_n) - \tO^{1,\alpha} (\U_n^{0}).\eea
Also set,
\bea
\label{0997}
  \tY_{N(n)} & := & Y_{N(n)}-\Exp\left[Y_{N(n)}|N(n)\right], \\
~~  \tY_{n-1-N(n)} & := & Y_{n-1-N(n)}-\Exp \left[Y_{n-1-N(n)}|N(n) \right] ,
\eea
and similarly for the $\tR$s and $\tS$s. Using the expressions for the expectations
at (\ref{0919wa}), (\ref{0204c}) and (\ref{0204cc}), from (\ref{0920iq}) we obtain
\bea
\label{0920jq}
\left( \begin{array}{l} \tY_n \\ \tR_n \\ \tS_n \end{array} \right)
 \eqd
\left( \begin{array}{lll} U^\alpha & 0 & 0 \\ 0 & 0 & 0 \\ 0 & U^\alpha & 0 \end{array} \right)
\left( \begin{array}{l} \tY^{\{1\}}_{N(n)} \\ \tR^{\{1\}}_{N(n)} \\
\tS^{\{1\}}_{N(n)} \end{array} \right)
\nonumber\\
+
\left( \begin{array}{lll} (1-U)^\alpha & 0 & 0 \\ 0 & (1-U)^\alpha & 0 \\ 0 & 0 & 0 \end{array} \right)
\left( \begin{array}{l} \tY^{\{2\}}_{n-1-N(n)} \\ \tR^{\{2\}}_{n-1-N(n)} \\
\tS^{\{2\}}_{n-1-N(n)} \end{array} \right) +
\left( \begin{array}{l} A_n \\ B_n \\ C_n
\end{array} \right),
\eea
where
\bean
 \left( \begin{array}{l} A_n \\ B_n \\ C_n \end{array} \right)
 =
\left( \begin{array}{l} \min \{ U, 1-U \}^\alpha + C_\alpha (n-1)^{(1/2)-\alpha} B_\alpha (n) +\frac{2^{-\alpha}}{\alpha-1} \left( U^\alpha +(1-U)^\alpha -1 \right)
\\ (U^\alpha-(1-U)^\alpha) \1_{\{ U > 1/2 \}} + \frac{1-2^{-\alpha}}{\alpha} ((1-U)^\alpha-1)
\\ (U^\alpha-\frac{1}{1+\alpha}) \frac{1-2^{-\alpha}-\alpha}{\alpha}
 \end{array} \right) \nonumber\\
+
\left( \begin{array}{l}
+U^\alpha h(N(n)) + (1-U)^\alpha h(n-1-N(n)) - h(n) \\
(1-U)^\alpha k (n-1-N(n)) -k(n) \\
U^\alpha k(N(n)) - \ell(n)
\end{array} \right)
,
\eean
where $B_\alpha(n)$ is as
defined
at (\ref{0803d}), $h(n)$, $k(n)$, $\ell(n)$ are all $o(1)$ as $n \to \infty$ and $C_\alpha$ is a constant.

In order to apply Theorem 4.1 of \cite{neinrusch},
we need to verify the conditions (24), (25) and (26) there.
By Lemma \ref{0803a}, $(n-1)^{(1/2)-\alpha}B_\alpha(n)$ tends to zero in $L^2$ as $n \to \infty$, for
$1/2<\alpha<1$. Thus, for condition (24) in \cite{neinrusch},
as $n \to \infty$,
\bea
\label{0920cbd}
\left( \begin{array}{l} A_n \\ B_n \\ C_n \end{array} \right) \inLL
\left( \begin{array}{l} \min \{ U, 1-U \}^\alpha  +\frac{2^{-\alpha}}{\alpha-1} \left( U^\alpha +(1-U)^\alpha -1 \right)
\\ (U^\alpha-(1-U)^\alpha)  \1_{\{ U > 1/2 \}} + \frac{1-2^{-\alpha}}{\alpha} ((1-U)^\alpha-1)
\\ (U^\alpha-\frac{1}{1+\alpha}) \frac{1-2^{-\alpha}-\alpha}{\alpha}
 \end{array} \right) .
\eea
Also, writing $\|\cdot\|_{\rm op}$
for the operator norm, for condition (25) in \cite{neinrusch},
\bea
\label{0920cad}
\Exp \left[ \left\| \left( \begin{array}{lll} U^\alpha & 0 & 0 \\ 0 & 0 & 0 \\ 0 & U^\alpha & 0 \end{array} \right)
\right\|_{\rm op}^2 + \left\| \left( \begin{array}{lll} (1-U)^\alpha & 0 & 0 \\ 0 & (1-U)^\alpha & 0 \\ 0 & 0 & 0 \end{array} \right)
\right\|_{\rm op}^2\right] = \frac{2}{2\alpha+1} <1,
\eea
for $\alpha>1/2$.
Finally, for condition (26) in \cite{neinrusch},
for $\alpha>0$ and any $\ell \in \N$, as $n \to \infty$
\bea
\label{0920ccd}
 \Exp \left[ \1_{ \{ N(n) \leq \ell \} \cup \{ N(n) = n \}}
U^{2\alpha} \right] \to 0; ~~
\Exp \left[ \1_{ \{ n-1-N(n) \leq \ell \} \cup \{ n-1-N(n) = n \}}
(1-U)^{2\alpha} \right] \to 0.
\eea
Taking $s=2$ and $C_n$ to be the identity matrix, Theorem
4.1 of \cite{neinrusch} applied to equation (\ref{0920jq}),
with the conditions
(\ref{0920cad}), (\ref{0920cbd}) and (\ref{0920ccd}),
implies that $(\tY_n,\tR_n,\tS_n)$ converges in Zolotarev
$\zeta_2$ metric (which implies convergence in
distribution; see e.g.~Chapter 14 of \cite{rachev})
 to $(\tY,\tR,\tS)$, where $\Exp[\tY]=\Exp[\tR]=\Exp[\tS]=0$
and the
distribution of $(\tY,\tR,\tS)$ is characterized by
the fixed-point equation
\bea
\label{0207av}
\left( \begin{array}{l} \tY \\ \tR  \\ \tS \end{array} \right)
& \eqd&
\left( \begin{array}{lll} U^\alpha & 0 & 0 \\ 0 & 0 & 0 \\ 0 & U^\alpha & 0 \end{array} \right)
\left( \begin{array}{l} \tY^{\{1\}} \\ \tR^{\{1\}}
\\ \tS^{\{1\}} \end{array} \right)
+
\left( \begin{array}{lll} (1-U)^\alpha & 0 & 0 \\ 0 & (1-U)^\alpha & 0 \\ 0 & 0 & 0 \end{array} \right)
\left( \begin{array}{l} \tY^{\{2\}} \\ \tR^{\{2\}} \\ \tS^{\{2\}}\end{array} \right) \nonumber\\ &  &
+
\left( \begin{array}{l} \min \{ U, 1-U \}^\alpha +((1-U)^\alpha +U^\alpha-1) \frac{2^{-\alpha}}{\alpha-1}
 \\ (U^\alpha-(1-U)^\alpha)  \1_{\{ U > 1/2 \}} + ((1-U)^\alpha-1) \frac{1-2^{-\alpha}}{\alpha} \\ \left(U^\alpha -\frac{1}{1+\alpha} \right) \frac{1-2^{-\alpha}-\alpha}{\alpha} \end{array}
 \right). \eea
 That is, $\tY$ satisfies
(\ref{0919z}), so that $\tY$ has the distribution of $\tJal$,
and setting $\tY=\tJal$ in (\ref{0207av}) gives (\ref{0207aav}). Then (\ref{0217qq})
follows by (\ref{0996}). $\square$

\subsection{Limit theory when $\alpha$=1}
\label{sec1}

Proposition \ref{0201a} below
is our main convergence result when $\alpha=1$.
First, we need the following
result,
the proof of which we defer to the Appendix. For $x \geq 0$, set
$\log^+ x := \max \{ \log x , 0\}$.

\begin{lemma}
\label{0703a}
Let $U$ be uniform on $(0,1)$ and, given $U$, let
$N(n) \sim {\rm Bin} (n-1,U)$. Then, as $n \to \infty$,
\bea
\label{0703c}
U (\log^+ N(n) - \log {n}) & \inLL & U \log U; \\
\label{0703d}
(1-U) (\log^+ (n-1-N(n))-\log {n}) & \inLL & (1-U) \log (1-U). \eea
\end{lemma}

\begin{proposition}
\label{0201a}
As $n \to \infty$,
\bea
\label{0217q}
\left( \begin{array}{l} \tO^{1,1} (\U_n^{0,1})
\\
\tO^{1,1} (\U_n^{0}) - \tO^{1,1} (\U_n^{0,1}) \\
\tO^{1,1} (\U_n) - \tO^{1,1} (\U_n^0)
\end{array}
\right)
\tod
\left( \begin{array}{l} \tJo
\\
\tR
\\
\tS
\end{array}
\right)
,\eea
where $(\tJo, \tR,\tS)$ satisfies the fixed-point equation
\bea
\label{0207aa}
\left( \begin{array}{l} \tJo \\ \tR  \\ \tS \end{array} \right)
& \eqd&
\left( \begin{array}{lll} U & 0 & 0 \\ 0 & 0 & 0 \\ 0 & U & 0 \end{array} \right)
\left( \begin{array}{l} \tJo^{\{1\}} \\ \tR^{\{1\}}
\\ \tS^{\{1\}} \end{array} \right)
+
\left( \begin{array}{lll} 1-U & 0 & 0 \\ 0 & 1-U & 0 \\ 0 & 0 & 0 \end{array} \right)
\left( \begin{array}{l} \tJo^{\{2\}} \\ \tR^{\{2\}} \\ \tS^{\{2\}}\end{array} \right) \nonumber\\ &  &
+
\left( \begin{array}{l} \frac{U}{2} \log U + \frac{1-U}{2} \log (1-U) +   \min \{ U, 1-U \}
 \\ (2U-1) \1_{\{U>1/2\}} - \frac{U}{2} \\ \frac{1}{4}- \frac{U}{2} \end{array} \right). \eea
In particular, $\tJo$ satisfies the fixed-point
equation (\ref{0919x}).
Also, $\Exp[ \tJo]=\Exp[\tR]=\Exp[\tS]=0$, $\Var[\tR]=1/16$, $\Var[\tS]=1/24$, and
\bea
\label{0920f}
\Var [ \tJo ] = \frac{1}{4} \left( 1 + \log {2} \right) - \frac{\pi^2}{24}
\approx 0.012053,
\eea
and $\Exp[ \tJo^3 ] \approx -0.00005732546$.
\end{proposition}
\proof We follow the proof of Proposition \ref{1003e}. Recall the definition
of $Y_n$ at (\ref{0999}). Again define $R_n$ and $S_n$ as at (\ref{0995}),
this time with $\alpha=1$.
Then we have that the $\alpha=1$ case of (\ref{0920iq}) holds.
We now renormalise (\ref{0920iq}), with the
notation of (\ref{0998}) and (\ref{0997}). By (\ref{0919wb}) we have
\[ \Exp[Y_n] = \Exp [\OO^{1,1} (\U_n^{0,1}) ] -1 = \frac{1}{2} \log n + \frac{1}{2}(\gamma-1) + h(n), \]
where $h(n)=o(1)$,
while by the $\alpha=1$ case of (\ref{0204c})
$\Exp[R_n] = (1/2) + k(n)$, where $k(n)=O(n^{-1})$,
and by $\alpha=1$ case of (\ref{0204cc}) $\Exp[S_n] = -(1/4) + \ell (n)$, where $\ell(n)=O(n^{-1})$.
Then by (\ref{0920iq})
\bea
\label{0920j}
\left(\!\! \begin{array}{l} \tY_n \\ \tR_n \\ \tS_n \end{array}\!\!\right)
 \eqd
\left(\!\! \begin{array}{lll} U & 0 & 0 \\ 0 & 0 & 0 \\ 0 & U & 0 \end{array} \!\!\right)
\!\!
\left(\!\! \begin{array}{l} \tY^{\{1\}}_{N(n)} \\ \tR^{\{1\}}_{N(n)} \\
\tS^{\{1\}}_{N(n)} \end{array} \!\!\right) \!
+ \!
\left( \!\!\begin{array}{lll} 1-U & 0 & 0 \\ 0 & 1-U & 0 \\ 0 & 0 & 0 \end{array} \!\!\right)
\!\!
\left( \!\!\begin{array}{l} \tY^{\{2\}}_{n-1-N(n)} \\ \tR^{\{2\}}_{n-1-N(n)} \\
\tS^{\{2\}}_{n-1-N(n)} \end{array} \!\!\right) \!
+\!
\left(\!\! \begin{array}{l} A_n \\ B_n \\ C_n
\end{array} \!\!\right),
\eea
where
\bean
 \left( \begin{array}{l} A_n \\ B_n \\ C_n \end{array} \right)
 =
\left( \begin{array}{l}
Uh(N(n)) +(1-U)h(n-1-N(n)) -h(n) \\
(1-U) k (n-1-N(n)) -k(n) \\
U k(N(n)) - \ell(n)
\end{array} \right)\\
+
\left( \begin{array}{l} \min \{ U, 1-U \} + \frac{U}{2} (\log^+ N(n)-\log n) +
\frac{1-U}{2} (\log^+ (n-1-N(n)) - \log n)
\\ (2U-1) \1_{\{ U > 1/2 \}} - \frac{U}{2}
\\ \frac{1}{4}-\frac{U}{2}
 \end{array} \right)
.
\eean
The conditions of Theorem 4.1 of \cite{neinrusch}
are satisfied, by (\ref{0920cbd}), (\ref{0920ccd}) and Lemma \ref{0703a}.
Taking $s=2$ and $C_n$ to be the identity, Theorem
4.1 of \cite{neinrusch} applied to equation (\ref{0920j})
shows that $(\tY_n,\tR_n,\tS_n)$ converges in Zolotarev
$\zeta_2$ metric and hence in
distribution
 to $(\tY,\tR,\tS)$, where $\Exp[\tY]=\Exp[\tR]=\Exp[\tS]=0$
and the
distribution of $(\tY,\tR,\tS)$ is characterized by
the fixed-point equation
\bea
\label{0207a}
\left( \begin{array}{l} \tY \\ \tR  \\ \tS \end{array} \right)
& \eqd&
\left( \begin{array}{lll} U & 0 & 0 \\ 0 & 0 & 0 \\ 0 & U & 0 \end{array} \right)
\left( \begin{array}{l} \tY^{\{1\}} \\ \tR^{\{1\}}
\\ \tS^{\{1\}} \end{array} \right)
+
\left( \begin{array}{lll} 1-U & 0 & 0 \\ 0 & 1-U & 0 \\ 0 & 0 & 0 \end{array} \right)
\left( \begin{array}{l} \tY^{\{2\}} \\ \tR^{\{2\}} \\ \tS^{\{2\}}\end{array} \right) \nonumber\\ &  &
+
\left( \begin{array}{l} \frac{U}{2} \log U + \frac{1-U}{2} \log (1-U) +   \min \{ U, 1-U \}
 \\ (2U-1) \1_{\{U>1/2\}} - \frac{U}{2} \\ \frac{1}{4}- \frac{U}{2} \end{array} \right).
\eea
That is, $\tY$ satisfies
(\ref{0919x}), so that $\tY$ has the distribution of $\tJo$,
and setting $\tY=\tJo$ in (\ref{0207a}) gives (\ref{0207aa}). By the $\alpha=1$ case of (\ref{0996})
we then have (\ref{0217q}).

It remains to prove the results for the higher moments of $\tJo$.
For the variance of $\tJo$, squaring both sides of
(\ref{0919x}), taking expectations, and using independence
and the
fact that $\Exp[ \tJo]=0$, we obtain
\bean \Exp [ \tJo^2 ] & = & \frac{2}{3} \Exp [ \tJo^2 ] + \Exp [ \min \{ U, 1-U \}^2 ]
+ \frac{1}{2} \Exp [ U^2 (\log U)^2 ]  \nonumber\\
& & + \frac{1}{2} \Exp[ U(1-U) \log U \log (1-U)]
+2 \Exp [ U \log U \min \{ U, 1-U \} ] .\eean
The integrals required for the expectations are standard,
and we find that $\Exp [ \tJo^2] = ((1+\log{2})/4) - (\pi^2/24)$,
which yields (\ref{0920f}). Similarly, we obtain the third moment
$\Exp[ \tJo^3] = -0.00005732546\ldots$
from (\ref{0919x}), although in this case numerical methods are required
for some of the integrals. $\square$

\subsection{Limit theory for $\alpha>1$}
\label{sec3}

\begin{proposition}
\label{0920c}
Let $\alpha>1$. \begin{itemize}
\item[(i)]
There exists a r.v.~$\Jal$
such that as $n \to \infty$
$\OO^{1,\alpha}(\U_n^{0,1}) \to 1+\Jal$ a.s.~and in $L^p$, $p \in \N$.
Also, $\Jal$
satisfies the fixed-point equality (\ref{0919y}), and
$\Exp[ \Jal] = 2^{-\alpha}/(\alpha-1)$.
\item[(ii)]
There exists a r.v.~$\Hal$
such that as $n \to \infty$
$\OO^{1,\alpha}(\U_n^0) \to \Hal$ a.s.~and in $L^2$.
Also, $\Hal$
satisfies the fixed-point equality (\ref{0923l}), and
$\Exp[ \Hal] = (1/\alpha)+2^{-\alpha}/(\alpha(\alpha-1))$.
\end{itemize}
\end{proposition}
\proof First we prove part (i).
Let $T_i$ be the length of the $i$th edge of the $\onng$ on $\U_n^{0,1}$,
  as defined at (\ref{0918b}).
Let $\Jal := \sum_{i=1}^\infty
T_i^\alpha$.
The sum converges almost surely since
it has non-negative terms
and, by (\ref{0919wc}), has finite expectation for $\alpha >1$. By a similar argument
as the Proof of Theorem \ref{onngthm} (ii) in Section \ref{prf1027}, the $L^p$ convergence follows
by H\"older's inequality and dominated convergence.

We now identify the limit.
We have (\ref{1029f}), this time for $\alpha>1$.
As $n \to \infty$, $N(n)$ and $n-N(n)$ both tend to
infinity almost surely, and so, by taking
$n \to \infty$ in (\ref{1029f}), we obtain
the fixed-point equation (\ref{0919y}).

The identity $E[\Jal] = 2^{-\alpha}(\alpha -1)^{-1}$
is obtained either from (\ref{0919wc}),
or by taking expectations in (\ref{0919y}).
 Next, if we set $\tJal = \Jal- \Exp[ \Jal]$, (\ref{0919y})
 yields (\ref{0919z}).

We now prove part (ii). Following the above argument
with the $H_i$ replacing the $T_i$ and using
(\ref{0923wc}) in place of (\ref{0919wc}) gives that
$\OO^{1,\alpha} (\U^0_n)$ converges a.s.~and in $L^p$, $p \in \N$,
to some random variable. Once more, we need to identify the limit.

Consider the $\alpha>1$ case of (\ref{0923m}).
As $n \to \infty$, $N(n)$ and $n-N(n)$ both tend to
infinity almost surely, and so, by taking
$n \to \infty$ in (\ref{0923m}), and using the
fact that $\OO^{1,\alpha} (\U_{N(n)}^{0,1})$
converges almost surely to $1+\Jal$ (by part (i)),
and that $\OO^{1,\alpha} (\U_{n-1-N(n)}^0)$
converges almost surely to $\Hal$ (by the argument above)
we obtain
the fixed-point equation (\ref{0923l}).

The identity $E[\Hal] = \alpha^{-1}+2^{-\alpha}\alpha^{-1}(\alpha -1)^{-1}$
is obtained either from (\ref{0923wc}),
or by taking expectations in (\ref{0923l}).
 Next, if we set $\tHal = \Hal- \Exp[ \Hal]$, (\ref{0923l})
 yields (\ref{0923n}).
$\square$

\subsection{Proof of Theorems \ref{onng01} and \ref{onng1}}
\label{sec4}

\noindent {\bf Proof of Theorem \ref{onng01}.} 
First we prove part (i) of the theorem.
For $1/2<\alpha<1$
we have that
\bea
\label{1006zz}
\left(\!\! \begin{array}{l} \tO^{1,\alpha}(\U_n^{0,1}) \\
\tO^{1,\alpha} (\U_n^0) \\
\tO^{1,\alpha} (\U_n) \end{array} \!\!\right)
= \left(\!\! \begin{array}{lll}
1 & 0 & 0 \\
1 & 1 & 0 \\
1 & 1 & 1 \end{array} \!\!\right)
\!\!
\left(\!\! \begin{array}{l} \tO^{1,\alpha}(\U_n^{0,1}) \\
\tO^{1,\alpha} (\U_n^0)-\tO^{1,\alpha}(\U_n^{0,1}) \\
\tO^{1,\alpha} (\U_n)-\tO^{1,\alpha}(\U_n^{0}) \end{array} \!\!\right)
\tod
\left(\!\! \begin{array}{lll}
1 & 0 & 0 \\
1 & 1 & 0 \\
1 & 1 & 1 \end{array} \!\!\right)
\!\!
\left(\!\! \begin{array}{l} \tJal \\
\tR \\
\tS \end{array} \!\!\right), \eea as $n \to \infty$, by Proposition
\ref{1003e}. By (\ref{0207aav}), the final term in (\ref{1006zz})
is equal in distribution to \bean \left(\!\! \begin{array}{lll}
1 & 0 & 0 \\
1 & 1 & 0 \\
1 & 1 & 1 \end{array} \!\!\right)
\!\!
\left(\!\! \begin{array}{lll}
U^\alpha & 0 & 0 \\
0 & 0 & 0 \\
0 & U^\alpha & 0 \end{array} \!\!\right)
\!\!
\left(\!\! \begin{array}{l} \tJal^{\{1\}} \\
\tR^{\{1\}} \\
\tS^{\{1\}} \end{array} \!\!\right) 
\!+\!
\left(\!\! \begin{array}{lll}
1 & 0 & 0 \\
1 & 1 & 0 \\
1 & 1 & 1 \end{array} \!\!\right)
\!\!
\left(\!\! \begin{array}{lll}
(1-U)^\alpha & 0 & 0 \\
0 & (1-U)^\alpha & 0 \\
0 & 0 & 0 \end{array} \!\!\right)
\!\!
\left(\!\! \begin{array}{l} \tJal^{\{2\}} \\
\tR^{\{2\}} \\
\tS^{\{2\}} \end{array} \!\!\right) \\
+ \left(\!\! \begin{array}{lll}
1 & 0 & 0 \\
1 & 1 & 0 \\
1 & 1 & 1 \end{array} \!\!\right)
\!\!
\left(
\!\! \begin{array}{l} \min \{ U, 1-U\}^\alpha +((1-U)^\alpha +U^\alpha-1) \frac{2^{-\alpha}}{\alpha-1} \\
(U^\alpha-(1-U)^\alpha) \1_{\{ U>1/2\}} +\frac{1-2^{-\alpha}}{\alpha} ((1-U)^\alpha-1) \\
(U^\alpha -\frac{1}{1+\alpha}) \frac{1-2^{-\alpha}-\alpha}{\alpha} \end{array} \!\!\right)
 .\eean
 Multiplying out and using the fact that $(U^\alpha-(1-U)^\alpha)\1_{\{U>1/2\}}=U^\alpha-\min\{U,1-U\}^\alpha$ we obtain
 \bean
\left( \begin{array}{l} \tO^{1,\alpha}(\U_n^{0,1}) \\
\tO^{1,\alpha} (\U_n^0) \\
\tO^{1,\alpha} (\U_n) \end{array} \right)
\tod
\left( \begin{array}{l} \tJal \\
 \tJal +\tR \\
\tJal+\tR+\tS \end{array} \right)  \eqd
U^\alpha \left( \begin{array}{l} \tJal^{\{1\}} \\
 \tJal^{\{1\}} \\
 \tJal^{\{1\}}+\tR^{\{1\}} \end{array} \right)
+ (1-U)^\alpha \left( \begin{array}{l}  \tJal^{\{2\}} \\
 \tJal^{\{2\}}+\tR^{\{2\}} \\
 \tJal^{\{2\}}+\tR^{\{2\}} \end{array} \right) \\
+ \left( \begin{array}{l} \min \{U,1-U\}^\alpha +((1-U)^\alpha +U^\alpha-1)\frac{2^{-\alpha}}{\alpha-1} \\
U^\alpha (1+\frac{2^{-\alpha}}{\alpha-1}) + ((1-U)^\alpha-1) ( \frac{1}{\alpha} + \frac{2^{-\alpha}}{\alpha(\alpha-1)}
)
 \\
(U^\alpha+(1-U)^\alpha-\frac{2}{1+\alpha})(\frac{1}{\alpha}-\frac{2^{-\alpha}}{\alpha(1-\alpha)}) \end{array} \right). \eean
So setting $\tHal=\tJal+\tR$ and
$\tGal=\tJal+\tR+\tS$, we have (\ref{1703b}).

Now we prove part (ii) of the theorem. For $\alpha=1$, as an analogoue of (\ref{1006zz}),
\bea
\label{1006z}
\left(\!\! \begin{array}{l} \tO^{1,1}(\U_n^{0,1}) \\
\tO^{1,1} (\U_n^0) \\
\tO^{1,1} (\U_n) \end{array} \!\!\right)
= \left(\! \begin{array}{lll}
1 & 0 & 0 \\
1 & 1 & 0 \\
1 & 1 & 1 \end{array} \!\right)
\!\!
\left(\!\! \begin{array}{l} \tO^{1,1}(\U_n^{0,1}) \\
\tO^{1,1} (\U_n^0)-\tO^{1,1}(\U_n^{0,1}) \\
\tO^{1,1} (\U_n)-\tO^{1,1}(\U_n^{0}) \end{array} \!\!\right)
\tod
\left(\! \begin{array}{lll}
1 & 0 & 0 \\
1 & 1 & 0 \\
1 & 1 & 1 \end{array} \!\right)
\!\!
\left(\!\! \begin{array}{l} \tJo \\
\tR \\
\tS \end{array} \!\!\right), \eea as $n \to \infty$, by Proposition
\ref{0201a}. By (\ref{0207aa}), the final term in (\ref{1006z}) is
equal in distribution to \bean \left( \!\begin{array}{lll}
1 & 0 & 0 \\
1 & 1 & 0 \\
1 & 1 & 1 \end{array}\! \right)
\!\!
\left(\! \begin{array}{lll}
U & 0 & 0 \\
0 & 0 & 0 \\
0 & U & 0 \end{array} \!\right)
\!\!
\left(\! \begin{array}{l} \tJo^{\{1\}} \\
\tR^{\{1\}} \\
\tS^{\{1\}} \end{array} \!\right) +
\left(\! \begin{array}{lll}
1 & 0 & 0 \\
1 & 1 & 0 \\
1 & 1 & 1 \end{array} \!\right)
\!
\left(\! \begin{array}{lll}
1-U & 0 & 0 \\
0 & 1-U & 0 \\
0 & 0 & 0 \end{array} \!\right)
\!
\left(\! \begin{array}{l} \tJo^{\{2\}} \\
\tR^{\{2\}} \\
\tS^{\{2\}} \end{array} \!\right) \\
+ \left(\! \begin{array}{lll}
1 & 0 & 0 \\
1 & 1 & 0 \\
1 & 1 & 1 \end{array}\!\right)
\!
\left(\! \begin{array}{l} \frac{U}{2} \log U + \frac{1-U}{2} \log (1-U)
+ \min \{ U, 1-U\} \\
(2U-1) \1_{\{ U>1/2\}} -\frac{U}{2} \\
\frac{1}{4}-\frac{U}{2} \end{array} \!\right)
 .\eean
 Multiplying out and using the fact that $(2U-1)\1_{\{U>1/2\}}=U-\min\{U,1-U\}$ we have
 \bean
\left( \begin{array}{l} \tO^{1,1}(\U_n^{0,1}) \\
\tO^{1,1} (\U_n^0) \\
\tO^{1,1} (\U_n) \end{array} \right)
\tod
\left( \begin{array}{l} \tJo \\
 \tJo +\tR \\
\tJo+\tR+\tS \end{array} \right)  \eqd
U \left( \begin{array}{l} \tJo^{\{1\}} \\
 \tJo^{\{1\}} \\
 \tJo^{\{1\}}+\tR^{\{1\}} \end{array} \right)
+ (1-U) \left( \begin{array}{l}  \tJo^{\{2\}} \\
 \tJo^{\{2\}}+\tR^{\{2\}} \\
 \tJo^{\{2\}}+\tR^{\{2\}} \end{array} \right) \\
 + \left( \begin{array}{l} (U/2) \log U + \frac{1-U}{2} \log (1-U) +\min \{U,1-U\} \\
(U/2) \log U + \frac{1-U}{2} \log (1-U) + \frac{U}{2} \\
(U/2) \log U + \frac{1-U}{2} \log (1-U) + \frac{1}{4} \end{array} \right). \eean
So setting $\tHo=\tJo+\tR$ and
$\tGo=\tJo+\tR+\tS$, we have (\ref{0202c}).
Proposition \ref{0201a} gives $\Exp[\tJo]=\Exp[\tR]=\Exp[\tS]=0$,
and so $\Exp[\tHo]=\Exp[\tGo]=0$ also. Proposition \ref{0201a} also gives
$\Var[\tJo]$. We obtain the higher moments of $\tHo$ and $\tGo$ from
(\ref{0924e})
and (\ref{1006b}). The stated covariances follow from the fixed point equation (\ref{0207aa})
and the moments given in Proposition \ref{0201a}.

Finally, part (iii) of the theorem is Proposition \ref{0920c}.
$\square$ \\

\noindent
{\bf Proof of Theorem \ref{onng1}.} Parts (i) and (ii) of the theorem follow
directly from the corresponding parts of Theorem \ref{onng01}. It remains to
prove part (iii) of the theorem. Suppose $\alpha>1$. Consider
the $\alpha>1$ case of (\ref{0926ii}).
We use
the fact that
$N(n)$ and $n-N(n)$ tend to infinity almost surely,
the independence given $U$ and $N(n)$, and the convergence in $L^p$ and almost
surely
of $\tO^{1,\alpha} (\U_n^0)$ (for $\alpha>1$) to obtain the result.
$\square$

\section{Proof of Theorem \ref{nng1d}}
\label{nngvar}

\noindent
{\bf Proof of Theorem \ref{nng1d}.}
We make use of the theory of Dirichlet spacings as discussed in Section \ref{spacings}.
Since the nearest-neighbour (directed)
graph joins each vertex (which sits at the endpoint of each spacing apart from the points 0 and 1) to
its nearest neighbour, we have, for $n \geq 3$
\bea
\label{1105a}
\LL_1^{1,\alpha}(\U_n)
 = (S^n_2)^\alpha+(S^n_n)^\alpha+
\sum_{i=2}^{n-1} \left( \min \{ S^n_i,S^n_{i+1} \} \right)^\alpha
.\eea
Now, from (\ref{1105a}), using exchangeability we have that
\bean
\Exp[ \LL_1^{1,\alpha} (\U_n) ] =
 2\Exp [ (S^n_1)^\alpha ] + (n-2) \Exp [ ( \min \{ S^n_1,S^n_2 \} )^\alpha ] ,
\eean
where, from (\ref{1105b}) and (\ref{0530g}) we have
\bea
\label{0516a}
\Exp \left[ \left( \min \{ S^n_1,S^n_2 \} \right)^\alpha \right]
& = & 2^{-\alpha} \Exp [ (S_1^n)^\alpha ] =
2^{-\alpha} \frac{\Gamma (\alpha+1) \Gamma (n+1)}{\Gamma (n+\alpha+1)}.\eea
Then (\ref{0530a}) follows. We now prove (\ref{0530b}).
Squaring both sides of (\ref{1105a}) and taking expectations, we have
\bean
& & \Exp \left[ \left( \LL_1^{1,\alpha} (\U_n) \right)^2 \right]  \\
\! & = & \!
 \sum_{i=2}^{n-1} \Exp \left[ \left( \min \{ S^n_i,S^n_{i+1} \} \right)^{2\alpha}
 \right]
+ 2 \sum_{i=3}^{n-1}
\sum_{j=2}^{i-1} \Exp \left[ \left( \min \{ S^n_i,S^n_{i+1} \} \right)^\alpha
\left( \min \{ S^n_j,S^n_{j+1} \} \right)^\alpha
\right] \\
\! &  & \! + \Exp [ (S_2^n)^{2\alpha} ] +  \Exp [
(S_{n}^n)^{2\alpha} ] +
2 \sum_{i=2}^{n-1} \Exp [ (S_2^n)^\alpha (\min \{ S_i^n, S_{i+1}^n \})^\alpha ] \\
\! &  & \!+ 2 \sum_{i=2}^{n-1} \Exp [ (S_{n}^n)^\alpha (\min \{
S_i^n, S_{i+1}^n \})^\alpha ] + 2 \Exp [ (S_2^n)^\alpha
(S_{n}^n)^\alpha ] .\eean
Then, by exchangeability, \bea
\label{0517a} \Exp \left[ \left( \LL_1^{1,\alpha} (\U_n) \right)^2 \right]  = 
 (n-2) \Exp \left[ \left( \min \{ S^n_1,S^n_2 \}
\right)^{2\alpha} \right] + 2 \Exp[ (S_1^n S_2^n )^\alpha] \nonumber\\
 + (n-3)(n-4) \Exp \left[ \left( \min \{ S^n_1,S^n_2 \}
\right)^{\alpha}
\left( \min \{ S^n_3,S^n_4 \} \right)^{\alpha}\right] \nonumber\\
 + 2(n-3) \Exp \left[ \left( \min \{ S^n_1,S^n_2 \} \right)^{\alpha}
\left( \min \{ S^n_2,S^n_3 \} \right)^{\alpha}\right] + 2 \Exp [ (S_1^n)^{2\alpha}]
\nonumber\\
 + 4(n-3) \Exp [ (S_1^n)^\alpha (\min \{ S_2^n , S_3^n \})^\alpha]
+4 \Exp [ (S_1^n)^\alpha (\min \{S_1^n,S_2^n\})^\alpha ] . \eea
Now, by (\ref{0516b}) and (\ref{0516c}) we have
\[ \Exp[ (S_1^n)^\alpha (\min \{S_2^n,S_3^n\})^\alpha ] =2^{-\alpha}
\frac{\Gamma (n+1) \Gamma (1+\alpha)^2}
{\Gamma (n+1+2\alpha)},\]
and, using (\ref{0516b}) this time with (\ref{0516d}) we obtain
\[ \Exp[ (\min \{S_1^n,S_2^n\})^\alpha (\min \{S_3^n,S_4^n\})^\alpha ] =2^{-2\alpha}
\frac{\Gamma (n+1) \Gamma (1+\alpha)^2}
{\Gamma (n+1+2\alpha)}.\]
Also we have that
\bean
 \Exp [ (S_1^n)^\alpha (\min \{S_1^n, S_2^n\})^\alpha ] & = &
\Exp [ (S_1^n)^{2\alpha} \1_{\{S_1^n<S_2^n\}}]
+ \Exp [ (S_1^n)^{\alpha} (S_2^n)^{\alpha} \1_{\{S_1^n>S_2^n\}}] \\
& = & \frac{1}{2} \Exp [ (\min \{ S_1^n, S_2^n \})^{2 \alpha} ] + \frac{1}{2}
\Exp [ (S_1^n)^\alpha (S_2^n)^{\alpha} ].\eean
Hence from (\ref{0516a}) and (\ref{0516b}) we obtain
\[ \Exp [ (S_1^n)^\alpha (\min \{S_1^n, S_2^n\})^\alpha ] =
\frac{1}{2} \left(2^{-2\alpha} \Gamma
(1+2\alpha)+\Gamma(1+\alpha)^2 \right) \frac{\Gamma(n+1)}{\Gamma(n+1+2\alpha)}.\]
The final term on the right hand side of (\ref{0517a}) that
we need to evaluate is
\bea
\label{0517b}
 \Exp [ (\min \{S_1^n, S_2^n\})^\alpha (\min \{S_2^n, S_3^n\})^\alpha ]
& = &
 \Exp [ (S_2^n)^{2\alpha}
 \1_{\{S_2^n < S_1^n, ~S_2^n < S_3^n \}}] \nonumber\\
& &
+ 4 \Exp [ (S_1^n)^\alpha (S_2^n)^\alpha
 \1_{\{S_1^n<S_2^n<S_3^n\}}].\eea
For the first term on the right of (\ref{0517b}), by (\ref{0516dd}) we have
\bean
 \Exp [ (S_2^n)^{2\alpha} \1_{\{S_2^n < S_1^n, ~S_2^n < S_3^n \}}]
& = & \frac{1}{3} \Exp [ (\min \{ S_1^n, S_2^n, S_3^n\})^{2\alpha}]\\
& = & 3^{-1-2\alpha} \frac{\Gamma (1+2\alpha) \Gamma(n+1)}{\Gamma (n+1+2\alpha)}.\eean
Now consider the second term on the right of (\ref{0517b}). By a direct computation using (\ref{0606d}), we have
\bean &  & \Exp [ (S_1^n)^\alpha (S_2^n)^\alpha
 \1_{\{S_1^n<S_2^n<S_3^n\}}] \\
& = &
n(n-1)(n-2)\int_0^{1/3} \ud y \int_y^{(1-y)/2} \ud x \int_x^{1-x-y} \ud z x^\alpha
y^\alpha (1-x-y-z)^{n-3} \\
& = &
n(n-1)\int_0^{1/3} \ud y \int_y^{(1-y)/2} x^\alpha y^\alpha (1-y-2x)^{n-2} \ud x,
\eean
which, via the change of variables $w=y+2x$ and Fubini's theorem is the same as
\[ n(n-1) 2^{-\alpha-1} \int_0^1 \ud w (1-w)^{n-2} \int_0^{w/3} y^\alpha (w-y)^\alpha \ud y.\]
Setting $t=3y/w$ reduces this to
\[ n(n-1) 6^{-\alpha-1} \int_0^1 w^{1+2\alpha} (1-w)^{n-2} \ud w \int_0^1 t^\alpha (1-(t/3))^\alpha \ud t.\]
Using (\ref{gint}) for the integral involving $w$, and the fact that (see, e.g., 15.3.1 in \cite{as})
for $a>0$,
\[ \int_0^1 t^{a-1} (1-(t/z))^{-b} \ud t = \frac{1}{a} \!~_2F_1(b,a;a+1;z) \]
for the integral involving $t$,
we obtain the expression for
 $J_{n,\alpha}$ as given by (\ref{0530c}). Then, by (\ref{0517a}) and the subsequent
calculations, we obtain (\ref{0530b}).
Finally, (\ref{0530d}) follows from (\ref{0530b}) by (\ref{stir}).
$\square$

\section*{Appendix: technical lemmas}

\noindent
{\bf Proof of Lemma \ref{0803a}.}
The result is trivial when $\alpha=1$ or $\alpha=0$.
Suppose $0<\alpha<1$. Suppose $n>1$. To ease notation, for the duration of this proof,
set $m=n-1$.
Then we have that for any $U \in (0,1)$ and
$0 \leq N(n) \leq m$,
\bea
\label{2003a}
 -1 \leq U^\alpha \left( \frac{N(n)}{m} \right)^{1-\alpha}
+ (1-U)^\alpha  \left( \frac{m-N(n)}{m} \right)^{1-\alpha} -1 \leq 0,\eea
so that in particular $|B_\alpha(n)| \leq n^{1/2}$ almost surely for $0 \leq \alpha \leq 1$.
Let \[ W_n := \frac{N(n)-mU}{\sqrt{ mU(1-U)}} ,\]
so that $\Exp[W_n]=0$, $\Exp[W_n^2]=1$, and
\[ \frac{N(n)}{mU} = 1 + W_n \sqrt{ \frac{1-U}{mU} }; ~~~ \frac{m-N(n)}{m (1-U)} = 1 - W_n 
\sqrt{ \frac{U}{m(1-U)}} .\]
 Then, by Taylor's theorem,
\bea
\label{0803c1}
 U^\alpha \left( \frac{N(n)}{m} \right)^{1-\alpha} & = &
U \left( 1 + (1-\alpha) W_n \sqrt {\frac {1-U}{mU}}
-R_1(n) W_n^2 \frac{1-U}{mU} \right) \\
\label{0803c2}
& = & U \left( 1 +R_2(n) W_n \sqrt {\frac {1-U}{mU}} \right) ,\eea
for remainder terms $R_1(n)$, $R_2(n)$ (which depend on $W_n$ and $U$).
Similarly, we have
\bea
\label{0803c3}
& &  (1-U)^\alpha \left( \frac{m-N(n)}{m} \right)^{1-\alpha} \nonumber\\
& = &
(1-U) \left( 1 -(1-\alpha) W_n \sqrt{\frac{U}{m(1-U)}} - R_3(n)
W_n^2 \frac{U}{m(1-U)} \right) \\
\label{0803c4}
& = & (1-U) \left( 1 -R_4(n) W_n \sqrt{\frac{U}{m(1-U)}} \right) .\eea
By the Lagrange form of the remainder in Taylor's theorem
and a continuity argument at $x=0$ there exists a constant $B \in (0,\infty)$ such that
for $\beta =1-\alpha$,
\[ 0 \geq \frac{(1+x)^\beta -1-\beta x}{x^2} \geq -B, \textrm{ and } 0 \leq \frac{(1+x)^\beta-1}{x} \leq B,\]
for all $x \geq -1$.
Thus we we have, for $i \in \{1,2,3,4\}$,
\bea
\label{0803b}
 0 \leq R_i(n) < C,\eea
for a finite positive constant $C$.

For $n>1$, $m=n-1$, let $E_n$ denote the event $m^{-3/4} < U < 1-m^{-3/4}$.
From (\ref{0803c1}) and (\ref{0803c3}) we obtain
\bean
\left| B_\alpha (n) \1_{E_n} \right|
 = 
\left|  -R_1(n) W_n^2 (1-U) m^{-1/2} -R_3(n) W_n^2 U m^{-1/2} \right|
\1_{D^3_n} \1_{E_n} 
 \leq C m^{-1/2} W_n^2 \1_{E_n},\eean
for some $C \in (0,\infty)$.
By a standard moment generating function calculation,
\bea
\label{0525a}
\Exp[ (N(n)-mU)^6 |U] & = & mU(1-U) \left[ 15 m^2 U^2(1-U)^2 -130mU^2 (1-U)^2
\right. \nonumber\\
&  & \left. +25mU(1-U) -30U (1-U) (1-2U)^2 +1 \right] \nonumber\\
& \leq & mU(1-U) (15 m^2 U^2(1-U)^2
+25mU(1-U) +1 ).\eea
By (\ref{0525a}) we have that
\[ \Exp [ W_n^6 \1_{E_n} ] \leq \Exp [(N(n)-mU)^6 m^{-3} U^{-3} (1-U)^{-3} |E_n ]
= O(1) ,\]
as $n \to \infty$, so
from (\ref{0803e}) we have that
\bea
\label{0803e}
 B_\alpha(n) \1_{E_n} \inLLL 0 .\eea
Also, from (\ref{0803c2}) and (\ref{0803c4}) we have,
\[
\left| B_\alpha (n) \1_{E_n^c} \right|
 =
\left| (R_2(n)-R_4(n)) W_n U^{1/2} (1-U)^{1/2} \right| \1_{E_n^c} 
, \]
and so using (\ref{0803b}) we have
\bea
\label{0525b}
 \left| B_\alpha (n) \1_{E_n^c} \right|
 \leq  C | W_n | U^{1/2} (1-U)^{1/2} \1_{E_n^c}
.\eea
Now, from (\ref{0525a}) we have that
\[ \Exp[ (W_n U^{1/2} (1-U)^{1/2} )^6 ]
=m^{-3} \Exp [ (N(n)-mU)^6]
= O(1) ,\]
as $n \to \infty$, so by Cauchy-Schwarz and the fact
that $\Pr[E_n^c] = O(n^{-3/4})$ we obtain from
(\ref{0525b}) that as $n\to \infty$
\bea
\label{2103j}
 \Exp \left[
\left| B_\alpha (n) \1_{E_n^c} \right|^3 \right]
\to 0 .\eea
So (\ref{0803e}) and (\ref{2103j}) complete the proof. $\square$\\

  \noindent
  {\bf Proof of Lemma \ref{0703a}.}
For $n \in \N$,
let $M_n :=\log^+ N(n) - \log {n}-\log U$. First, suppose $N(n) \geq nU/2$.
We have that
\[  -\log 2 \leq M_n \1_{\{N(n) \geq nU/2\}} \1_{\{nU \geq 2\}} \leq -\log U .\]
Hence
\bea
\label{0703b}
 U^2 M_n^2 \1_{\{N(n) \geq nU/2\}} \1_{\{nU \geq 2\}} \leq U^2 \max \{ (\log 2)^2, (\log U)^2 \}.\eea
The expected value of the right hand side of (\ref{0703b}) is finite.
Also, $U^2 M_n^2 \toas 0$ as $n \to \infty$, by continuity and the strong law of
large numbers for $N(n)$. Hence, by the dominated convergence theorem,
\bea
\label{0703e}
 \Exp [ U^2 M_n^2 \1_{\{N(n) \geq nU/2\}} \1_{\{nU \geq 2\}} ] \to 0.\eea
Also, we have $0 \leq \log^+ N(n) \leq \log {n}$,
so that $-\log n \leq M_n \leq -\log U$.
Hence
\bea
\label{0905a}
 U^4 M_n^4 \leq (\log n )^4 + (\log U)^4 ,\eea
so that $\Exp [U^4 M_n^4] = O( (\log n)^4)$.
Since $\Pr[ nU < 2] = 2n^{-1}$, we then obtain, by Cauchy-Schwarz, that
there exists a finite positive constant $C$ such that
\bea
\label{0703h}
 \Exp [ U^2 M_n^2 \1_{\{N(n) \geq nU/2\}} \1_{\{nU < 2\}} ]
 \leq C (\log n)^2 n^{-1/2} \to 0,\eea
as $n \to \infty$. Now, suppose $0 \leq N(n) < nU/2$. In this case,
from (\ref{0905a}), and Cauchy-Schwarz again, for some finite positive constant $C$
\bea
\label{0703f}
 \Exp [ U^2 M_n^2 \1_{\{N(n) < nU/2\}} ] \leq C (\log {n})^2 (\Pr[N(n) < nU/2])^{1/2} \to 0,
\eea
as $n\to \infty$, since
\[ (\log n)^2 (\Pr[N(n) \! < nU/2])^{1/2}\! \leq
(\log n)^2  (  
 \Pr[ U \! < n^{-1/2}] +
 \Pr[ U \! >n^{-1/2}\!,  N(n) \! < nU/2] )^{1/2} \!\! ,\]
 which tends to zero as $n \to \infty$, using standard
bounds for the tail of a binomial distribution (see, e.g., Lemma 1.1 in \cite{penbook})
for the final probability.
The results (\ref{0703e}), (\ref{0703h}), and (\ref{0703f})
then give (\ref{0703c}).
The argument for (\ref{0703d}) is similar. $\square$ \\

\begin{center} \textbf{Acknowledgements} \end{center}
AW began this work while at the University of Durham, supported by an EPSRC
doctoral training account.

\end{document}